\documentclass[a4paper,12pt,twoside]{article}
\usepackage{amsmath,amssymb,ifthen}
\usepackage[amsmath,thmmarks]{ntheorem}
\usepackage{paper}
\usepackage[all]{xy}
\usepackage{mathrsfs}
\usepackage{url}

\newcommand{\red}{\mathrm{red}}

\DeclareMathOperator{\Spf}{Spf}

\DeclareMathOperator{\Spa}{Spa}
\DeclareMathOperator{\Nilp}{\mathbf{Nilp}}
\DeclareMathOperator{\Irr}{Irr}
\DeclareMathOperator{\cInd}{c-Ind}
\DeclareMathOperator{\height}{height}
\DeclareMathOperator{\slope}{slope}
\DeclareMathOperator{\simil}{sim}
\DeclareMathOperator{\Lie}{Lie}
\DeclareMathOperator{\ord}{ord}
\DeclareMathOperator{\Ig}{Ig}
\DeclareMathOperator{\GSp}{GSp}
\DeclareMathOperator{\GU}{GU}
\DeclareMathOperator{\inv}{inv}
\DeclareMathOperator{\Trd}{Trd}
\newcommand{\M}{\mathcal{M}}
\newcommand{\X}{\mathbb{X}}

\SelectTips{cm}{11}
\title{On irreducible components of Rapoport-Zink spaces}
\author{Yoichi Mieda}

\begin{document}
\maketitle

\begin{firstfootnote}
 Graduate School of Mathematical Sciences, The University of Tokyo, 3--8--1 Komaba, Meguro-ku, Tokyo, 153--8914, Japan

 E-mail address: \texttt{mieda@ms.u-tokyo.ac.jp}

\end{firstfootnote}

\begin{abstract}
 Under a mild condition, we prove that the action of the group of self-quasi-isogenies
 on the set of irreducible components of a Rapoport-Zink space has finite orbits. 
 Our method allows both ramified and non-basic cases.
 As a consequence, we obtain some finiteness results on the representation obtained from the $\ell$-adic
 cohomology of a Rapoport-Zink tower.
\end{abstract}

\section{Introduction}
A Rapoport-Zink space, introduced in \cite{MR1393439}, is a moduli space of deformations
by quasi-isogenies of a fixed $p$-divisible group $\X$ over $\overline{\F}_p$
with various additional structures.
It is a formal scheme $\M$ locally formally of finite type over the ring of integers
$\mathcal{O}_{\breve{E}}$ of $\breve{E}$, which is the completion of the maximal unramified extension of
a finite extension $E$ (called the local reflex field) of $\Q_p$.
Rapoport-Zink spaces are local analogues of Shimura varieties
of PEL type (\cf \cite{MR1124982}), and are also related to Shimura varieties themselves
by the theory of $p$-adic uniformization (\cf \cite[Chapter 6]{MR1393439}).

One of the main reasons we are interested in Rapoport-Zink spaces is the conjectural relation
between the $\ell$-adic cohomology of the Rapoport-Zink tower and the local Langlands correspondence
(\cf \cite[Conjecture 5.1]{MR1403942}).
Let us recall it briefly. We write $M$ for the rigid generic fiber of $\M$.
By using level structures on the universal $p$-divisible group over $M$, we can construct
a projective system of \'etale coverings $\{M_{K'}\}$ over $M$ called the Rapoport-Zink tower.
Here $K'$ runs through compact open subgroups of $G'=\mathbf{G}'(\Q_p)$,
where $\mathbf{G}'$ is a certain inner form of the reductive algebraic group $\mathbf{G}$ over $\Q_p$
that is naturally attached to the linear-algebraic data appearing in the definition of $\M$.
As in the case of Shimura varieties, $G'$ acts on the tower by Hecke correspondences.
On the other hand, the group $J$ of self-quasi-isogenies of $\X$ preserving
additional structures naturally acts on $\M$, and this action lifts canonically to $M_{K'}$
for each $K'\subset G'$.
Let $H^i_c(M_\infty)$ be the compactly supported $\ell$-adic \'etale cohomology of the tower
$\{M_{K'}\}$.
This is naturally equipped with an action of $G'\times J\times W_E$, where $W_E$ is the Weil
group of $E$. Roughly speaking, it is expected that the representation
$H^i_c(M_\infty)$ of $G'\times J\times W_E$ is described by means of
the local Langlands correspondence for $G'$ and $J$.
There are many results in the classical setting, namely when $\M$ is either the Lubin-Tate space
or the Drinfeld upper half space (\cf \cite{MR1044827}, \cite{MR1464867}, \cite{MR1876802},
\cite{MR2511742}, \cite{MR2308851}). However, very little is known in other cases.

In this paper, we will study the underlying reduced scheme $\overline{\M}=\M^{\red}$.
There already have been a lot of results on $\overline{\M}$ in various particular cases;
see \cite{MR2369090}, \cite{MR2466182}, \cite{MR2800696}, \cite{MR3175176}
and \cite{MR3272278} for example.
In contrast to them, here we try to keep $\M$ as general as possible.
The main theorem of this paper is as follows.

\begin{thm}[Theorem \ref{thm:main}]\label{thm:main-intro}
 Assume that the isogeny class of the $p$-divisible group $\X$
 with additional structures comes from an abelian variety (for the precise definition, see
 Definition \ref{defn:comes-from-AV}). Then, the action of $J$ on the set of irreducible
 components of $\overline{\M}$ has finite orbits.
\end{thm}

As a consequence of this theorem, we can obtain the following finiteness result on
the cohomology $H^i_c(M_\infty)$.

\begin{thm}[Theorem \ref{thm:finiteness}]\label{thm:coh-finiteness}
 Assume the same condition as in Theorem \ref{thm:main-intro}.
 Then, for every compact open subgroup $K'$ of $G'$, the $K'$-invariant part $H^i_c(M_\infty)^{K'}$ of
 $H^i_c(M_\infty)$ is a finitely generated $J$-representation.
\end{thm}

By combining Theorem \ref{thm:coh-finiteness} with the duality theorem proved in 
\cite{2017arXiv170906651K}, \cite{2017arXiv171006935C} and \cite{Scholze-Berkeley}, we obtain
the following finiteness result for the basic Rapoport-Zink tower for $\GSp_{2n}$.

\begin{thm}[Theorem \ref{thm:coh-fin-length}]\label{thm:coh-fin-length-intro}
 Assume that $\mathcal{M}$ is the basic Rapoport-Zink space for $\GSp_{2n}$,
 in which case $G'=G=\GSp_{2n}(\Q_p)$ and $J$ is a quaternion unitary similitude group.
 For every integers $i,r\ge 0$ and every irreducible smooth representation $\rho$ of $J$,
 the $G$-representation
 $\Ext^r_J(H^i_c(M_\infty),\rho)^{\text{$\mathcal{D}_c$-sm}}$ has finite length.
 In particular, the $G'$-representation
 $\Hom_J(H^i_c(M_\infty),\rho)^{\mathrm{sm}}$ has finite length.
\end{thm}
For the precise notation, see Section \ref{subsec:GSp}.
To prove this theorem, we use Theorem \ref{thm:coh-finiteness} for two Rapoport-Zink towers:
one is the basic Rapoport-Zink tower for $\GSp_{2n}$,
and the other is a Rapoport-Zink tower for a quaternion unitary similitude group.
These two towers are isomorphic at the infinite level, thanks to the duality isomorphism.

These finiteness results are very useful when we apply representation theory to the study of $H^i_c(M_\infty)$.
For example, in \cite[\S 2--3]{LJLC}, Theorem \ref{thm:coh-finiteness} for the Drinfeld tower, which is well-known,
enables the author to apply a result in \cite{MR1779801} to investigate $H^i_c(M_\infty)$
by using the Lefschetz trace formula. Further, by using Theorem \ref{thm:coh-fin-length-intro}
for the Drinfeld tower (\cite[Corollary 4.3]{LJLC}),
we obtained a purely local and geometric proof of the local Jacquet-Langlands correspondence in some cases
(see \cite[Theorem 6.10]{LJLC}).
In \cite{RZ-LTF}, the author extends these arguments to the basic Rapoport-Zink space for $\GSp_4$.
For this purpose, Theorem \ref{thm:coh-fin-length-intro} plays a crucial role.
We also remark that Theorem \ref{thm:main-intro} and Theorem \ref{thm:coh-fin-length-intro} are
indispensable in the study of the relation between the Zelevinsky involution and the cohomology
$H^i_c(M_\infty)$ (see \cite{RZ-Zel});
in fact, Theorem \ref{thm:main-intro} is the main condition for applying the general theory developed in \cite{RZ-Zel},
and Theorem \ref{thm:coh-fin-length-intro} is needed because
the duality theory in \cite{MR1471867} behaves well only for representations of finite length.

Note that the same results as Theorem \ref{thm:main-intro} and Theorem \ref{thm:coh-finiteness}
have been obtained by Fargues (\cite[Th\'eor\`eme 2.4.13, Proposition 4.4.13]{MR2074714})
under the condition that the data defining $\M$ is unramified. 
However, there are many interesting ramified settings, such as the quaternion unitary case appearing in the proof of
Theorem \ref{thm:coh-fin-length-intro}.
Usually the geometry of $\M$ attached to a ramified data is much more complicated than that in the unramified case.
For example, in the unramified case $\M$ is formally smooth over $\mathcal{O}_{\breve{E}}$,
while in the ramified case $\M$ is sometimes not even flat over $\mathcal{O}_{\breve{E}}$.

We would like to say a few words on the condition on $\X$ in Theorem \ref{thm:main-intro}.
First, the problem to find an abelian variety with additional structures
which has the prescribed $p$-divisible group up to isogeny
can be seen as a generalization of Manin's problem, and is studied by many people
mainly under the unramified setting (see Remark \ref{rem:known-cases}).
As an example of ramified cases, we will prove that
any $\mathbb{X}$ comes from an abelian variety when $\mathbf{G}$ is a quaternion unitary similitude group;
see Proposition \ref{prop:GU(n,D)-comes-from-AV}.
By these results, it is natural for the author to expect that this condition always holds,
at least if the isogeny class of $\X$ is admissible in the sense of \cite[Definition 1.18]{MR1393439}.
Second, this condition is obviously satisfied if $\M$ is related to a Shimura variety
of PEL type by the theory of $p$-adic uniformization. At present, the most successful way to
study $H^i_c(M_\infty)$ is relating it with the cohomology of a suitable Shimura variety
and applying a global automorphic method. Rapoport-Zink spaces to which this technique is applicable
automatically satisfy our condition.
The author hopes that these two points guarantee
the condition in Theorem \ref{thm:main-intro} to be harmless in practice.

To explain the strategy of our proof of Theorem \ref{thm:main-intro}, first let us
assume that the isocrystal with additional structures $b$ associated with $\X$ is basic
(\cf \cite[\S 5.1]{MR809866}). In this case, by the $p$-adic uniformization theorem of Rapoport-Zink
\cite[Theorem 6.30]{MR1393439}, $\M$ uniformizes some moduli space of PEL type $X$,
whose generic fiber is often a disjoint union of Shimura varieties,
along an open and closed subscheme $Z$ of the basic stratum $\overline{X}{}^{(b)}$
of the special fiber $\overline{X}=X\otimes\overline{\F}_p$.
More precisely, there exist finitely many torsion-free discrete cocompact subgroups 
$\Gamma_1,\ldots,\Gamma_m$
of $J$ such that the formal completion $X_{/Z}$ of $X$ along $Z$ is isomorphic to
$\coprod_{i=1}^m\Gamma_i\backslash\M$. Since $\overline{X}$ is a scheme of finite type
over $\overline{\F}_p$, $Z$ has only finitely many irreducible components.
Hence the number of irreducible components of $\Gamma_1\backslash\overline{\M}$ is also finite,
and thus in particular the number of $J$-orbits in the set of irreducible components of $\overline{\M}$
is finite.
However, if $b$ is not basic, the $p$-adic uniformization theorem \cite[Theorem 6.23]{MR1393439}
becomes more complicated; it involves the formal completion of $X$ along a possibly infinite set
of closed subschemes of $\overline{X}$. For this reason, we cannot extend the method above to
the non-basic case. To overcome this problem, we will invoke results
by Oort \cite{MR2051612} and Mantovan \cite{MR2074715}, \cite{MR2169874}.
Roughly speaking, they proved that the Newton polygon stratum $\overline{X}{}^{(b)}$ is
almost equal to the product of $\overline{\M}$ with an appropriate Igusa variety (\cf \cite[\S 4]{MR2169874}).
We will generalize it to the ramified case.
Combining this generalization with the fact that $\overline{X}{}^{(b)}$ is of finite type,
we can conclude Theorem \ref{thm:main-intro}.
The author thinks that this method itself has some importance, because it is a first step
to a generalization of Mantovan's formula (the main theorem of \cite{MR2169874}) to the ramified case.
Recall that Mantovan's formula is a natural extension of the method in \cite{MR1876802},
and is one of the most powerful tools to investigate the $\ell$-adic cohomology of
Shimura varieties. It is also useful to study the cohomology $H^i_c(M_\infty)$ of Rapoport-Zink
towers (\cf \cite{MR2905002}).
We plan to pursue this problem in our future work.

The outline of this paper is as follows.
In Section \ref{sec:RZ-space}, we recall the definition of Rapoport-Zink spaces
and explain the precise statement of the main theorem.
We also observe that the condition on $\mathbb{X}$ in Theorem \ref{thm:main-intro} holds
in the quaternion unitary case.
Section \ref{sec:proof-main-thm} is devoted to the proof of the main theorem.
In Section \ref{sec:applications}, we will give some applications of the main theorem,
including Theorem \ref{thm:coh-finiteness} and Theorem \ref{thm:coh-fin-length-intro}.

\medbreak
\noindent{\bfseries Notation}\quad
For a field $k$, we denote its algebraic closure by $\overline{k}$.
For a scheme $X$ and a point $x$ in $X$, we write $\kappa(x)$ for the residue field at $x$.

\section{Rapoport-Zink spaces}\label{sec:RZ-space}
\subsection{Statement of the main theorem}\label{subsec:statement}
Here we briefly recall the definition of Rapoport-Zink spaces.
The main reference is \cite[Chapter 3]{MR1393439}.
Fix a prime number $p>2$.
A Rapoport-Zink space of PEL type is associated with a tuple $(B,\mathcal{O}_B,*,V,\langle\ ,\ \rangle,\mathscr{L},b,\mu)$
consisting of the following objects (\cf \cite[Definition 3.18]{MR1393439}):
\begin{itemize}
 \item $B$ is a finite-dimensional semisimple algebra over $\Q_p$.
 \item $\mathcal{O}_B$ is a maximal order of $B$.
 \item $*$ is an involution on $B$ under which $\mathcal{O}_B$ is stable.
 \item $V$ is a finite faithful $B$-module.
 \item $\langle\ ,\ \rangle\colon V\times V\longrightarrow \Q_p$ is a non-degenerate alternating bilinear pairing
       satisfying $\langle av,w\rangle=\langle v,a^*w\rangle$ for every $a\in B$ and $v,w\in V$.
 \item $\mathscr{L}$ is a self-dual multi-chain of $\mathcal{O}_B$-lattices in $V$
       (\cf \cite[Definition 3.1, Definition 3.13]{MR1393439}).
\end{itemize}
To explain the remaining objects, we denote by $\mathbf{G}$ the algebraic group over $\Q_p$
consisting of $B$-linear automorphisms of $V$
which preserve $\langle\ ,\ \rangle$ up to a scalar multiple (\cf \cite[1.38]{MR1393439}).
\begin{itemize}
 \item $b$ is an element of $\mathbf{G}(K_0)$, where $K_0$ denotes the fraction field
       of $W(\overline{\F}_p)$.
 \item $\mu$ is a cocharacter $\mathbb{G}_m\longrightarrow \mathbf{G}$
       defined over a finite extension $K$ of $K_0$.
\end{itemize}
We impose the following conditions:
\begin{enumerate}
 \item[(a)] The isocrystal $(N_b,\Phi_b)=(V\otimes_{\Q_p}K_0,b\sigma)$ has slopes in the interval $[0,1]$,
	   where $\sigma$ is the Frobenius isomorphism on $K_0$.
 \item[(b)] We have $\simil(b)=p$, where $\simil\colon \mathbf{G}\longrightarrow \mathbb{G}_m$ denotes the similitude character.
 \item[(c)] The weight decomposition of $V\otimes_{\Q_p}K$ with respect to $\mu$ has only the weight $0$ and $1$ parts:
       $V\otimes_{\Q_p}K=V_0\oplus V_1$.
\end{enumerate}
We call such a tuple a Rapoport-Zink datum.

Before proceeding, it is convenient to introduce notion of polarized $B$-isocrystals
(\cf \cite[Example 4.9]{MR2484281}).

\begin{defn}\label{defn:pol-isoc}
 A polarized $B$-isocrystal is a quadruple $(N,\Phi,\langle\ ,\ \rangle,c)$ consisting of
 \begin{itemize}
  \item an isocrystal $(N,\Phi)$ endowed with a $B$-action,
  \item a non-degenerate alternating bilinear pairing
	$\langle\ ,\ \rangle\colon N\times N\longrightarrow K_0$,
  \item and an element $c\in K_0^\times$
 \end{itemize}
 such that $\langle ax,y\rangle=\langle x,a^*y\rangle$ and
 $\langle \Phi x,\Phi y\rangle=c\sigma(\langle x,y\rangle)$ for each $a\in B$ and $x,y\in N$.
 We also refer to $(N,\Phi,\langle\ ,\ \rangle,c)$ as a $c$-polarized $B$-isocrystal $(N,\Phi,\langle\ ,\ \rangle)$.
 A morphism between two polarized $B$-isocrystals $(N,\Phi,\langle\ ,\ \rangle,c)$ and
 $(N',\Phi',\langle\ ,\ \rangle',c')$ is a morphism $f\colon N\longrightarrow N'$ of isocrystals
 compatible with $B$-action 
 such that there exists $a\in K_0^\times$ satisfying
 $\langle f(x),f(y)\rangle'=a\langle x,y\rangle$ for every $x,y\in N$.
 If $N\neq 0$, such an element $a$ satisfies $ac=c'\sigma(a)$.
 In particular, if $c=c'$, then $a$ lies in $\Q_p^\times$.

 Consider a triple $(X,\lambda,\iota)$ consisting of a $p$-divisible group $X$ over $\overline{\F}_p$,
 a quasi-polarization $\lambda\colon X\longrightarrow X^\vee$
 (namely, a quasi-isogeny satisfying $\lambda^\vee=-\lambda$)
 and a homomorphism $\iota\colon B\longrightarrow \End(X)\otimes_{\Z_p}\Q_p$
 satisfying $\lambda\circ\iota(a^*)=\iota(a)^\vee\circ\lambda$ for every $a\in B$.
 Then the rational Dieudonn\'e module $\mathbb{D}(X)_{\Q}$ is naturally endowed with a structure of 
 a $p$-polarized $B$-isocrystal.
 We denote it by $\mathbb{D}(X,\lambda,\iota)_{\Q}$, or simply by $\mathbb{D}(X)_{\Q}$.
\end{defn}

An element $b\in \mathbf{G}(K_0)$ gives rise to a polarized $B$-isocrystal
\[
 (N_b,\Phi_b,\langle\ ,\ \rangle,\simil(b)).
\]
We call it the polarized $B$-isocrystal attached to $b$ and simply write $N_b$ for it.
Let $\mathbf{J}$ be the algebraic group over $\Q_p$ of automorphisms of $N_b$;
for a $\Q_p$-algebra $R$, $\mathbf{J}(R)$ consists of $h\in \Aut_R(N_b\otimes_{\Q_p}R)$ such that
\begin{itemize}
 \item $h\circ (\Phi_b\otimes \id)=(\Phi_b\otimes \id)\circ h$,
 \item and $h$ preserves the pairing $\langle\ ,\ \rangle$ on $N_b\otimes_{\Q_p}R$ up to $R^\times$-multiplication.
\end{itemize}
For the representability of $\mathbf{J}$, see \cite[Proposition 1.12]{MR1393439}.
Moreover, if we denote by $\nu\colon \mathbb{D}\longrightarrow \mathbf{G}\otimes_{\Q_p}K_0$
the slope homomorphism attached to $b$ (\cf \cite[\S 4.2]{MR809866}) and by $\mathbf{G}_\nu$ the centralizer
of the image of $\nu$ in $\mathbf{G}\otimes_{\Q_p}K_0$,
then we have an isomorphism
$\mathbf{J}\otimes_{\Q_p}\overline{K}_0\yrightarrow{\cong}\mathbf{G}_\nu\otimes_{K_0}\overline{K}_0$
(see \cite[\S 3.3, \S A.2]{MR1485921}). In particular, $\mathbf{J}$ is reductive.
On the other hand, $\mathbf{J}$ is not connected in general.

The conditions (a) and (b) ensure the existence of a $p$-divisible group $\X$ over $\overline{\F}_p$
whose rational Dieudonn\'e module $\mathbb{D}(\X)_{\Q}$ is isomorphic to $(N_b,\Phi_b)$.
We fix such an $\X$.
Then, the pairing $\langle\ ,\ \rangle\colon N_b\times N_b\longrightarrow K_0$ induces an isomorphism of
isocrystals $\mathbb{D}(\X)_\Q\yrightarrow{\cong} \mathbb{D}(\X^\vee)_\Q$.
As the Dieudonn\'e functor $\mathbb{D}(-)$ is fully faithful, we obtain a quasi-isogeny $\lambda_0\colon \X\longrightarrow \X^\vee$, which is easily checked to be a quasi-polarization.
Similarly, we have a homomorphism
$\iota_0\colon B\longrightarrow \End((N_b,\Phi_b))\cong \End(\mathbb{D}(\X)_\Q)\cong \End(\X)\otimes_{\Z_p}\Q_p$,
which satisfies $\lambda\circ\iota(a^*)=\iota(a)^\vee\circ\lambda$ for every $a\in B$.
In this way, we find a triple $(\X,\lambda_0,\iota_0)$ as in Definition \ref{defn:pol-isoc} such that
$\mathbb{D}(\X,\lambda_0,\iota_0)_{\Q}$ is isomorphic to $N_b$ as a polarized $B$-isocrystal.
Again by the fully faithfulness of the Dieudonn\'e functor, the group $J=\mathbf{J}(\Q_p)$ can be identified with the group of self-quasi-isogenies of $\X$
compatible with $\iota_0$ and preserving $\lambda_0$ up to $\Q_p^\times$-multiplication.

Let $E$ be the field of definition of the cocharacter $\mu$. It is the subfield of $K$ generated by
$\{\Tr(a;V_0)\mid a\in B\}$ over $\Q_p$. In particular it is a finite extension of $\Q_p$. We denote by $\breve{E}$ the composite field
of $E$ and $K_0$ inside $K$. We write $\Nilp_{\mathcal{O}_{\breve{E}}}$ for the category of $\mathcal{O}_{\breve{E}}$-schemes
on which $p$ is locally nilpotent. For an object $S$ in $\Nilp_{\mathcal{O}_{\breve{E}}}$, we put $\overline{S}=S\otimes_{\mathcal{O}_{\breve{E}}}\mathcal{O}_{\breve{E}}/p\mathcal{O}_{\breve{E}}$.
We denote the category of sets by $\mathbf{Set}$.

Now we can give the definition of the Rapoport-Zink space (see also \cite[Definition 3.21, 3.23 c), d)]{MR1393439}):

\begin{defn}\label{defn:RZ-space}
 Consider the functor $\M\colon \Nilp_{\mathcal{O}_{\breve{E}}}\longrightarrow \mathbf{Set}$ that associates $S$ with
 the set of isomorphism classes of $\{(X_L,\iota_L,\rho_L)\}_{L\in \mathscr{L}}$ where
 \begin{itemize}
  \item $X_L$ is a $p$-divisible group over $S$,
  \item $\iota_L\colon \mathcal{O}_B\longrightarrow \End(X_L)$ is a homomorphism,
  \item and $\rho_L\colon \X\otimes_{\overline{\F}_p}\overline{S}\longrightarrow X_L\times_S\overline{S}$
	is a quasi-isogeny compatible with $\mathcal{O}_B$-actions,
 \end{itemize}
 such that the following conditions are satisfied.
\begin{enumerate}
 \item[(a)] For $L,L'\in\mathscr{L}$ with $L\subset L'$, the quasi-isogeny $\rho_{L'}\circ \rho^{-1}_L\colon X_L\times_S\overline{S}\longrightarrow X_{L'}\times_S\overline{S}$ lifts to an isogeny $X_L\longrightarrow X_{L'}$ with height $\log_p\#(L'/L)$.
	    Such a lift is automatically unique and $\mathcal{O}_B$-equivariant.
	    We denote it by $\widetilde{\rho}_{L',L}$.
 \item[(b)] For $a\in B^\times$ which normalizes $\mathcal{O}_B$ and $L\in\mathscr{L}$, the quasi-isogeny
	    $\rho_{aL}\circ \iota_0(a)\circ \rho_L^{-1}\colon X^a_L\times_S\overline{S}\longrightarrow X_{aL}\times_S\overline{S}$
	    lifts to an isomorphism $X^a_L\longrightarrow X_{aL}$. Here $X^a_L$ denotes the $p$-divisible group $X_L$
	    with the $\mathcal{O}_B$-action given by $\iota^a_L\colon \mathcal{O}_B\longrightarrow \End(X_L)$; $x\longmapsto \iota_L(a^{-1}xa)$.
 \item[(c)] Locally on $S$, there exists a constant $c\in\Q_p^\times$ such that for every $L\in\mathscr{L}$, we can find an isomorphism
	    $p_L\colon X_L\longrightarrow (X_{L^\vee})^\vee$ which makes the following diagram commute:
	    \[
	     \xymatrix{%
	    \X\otimes_{\overline{\F}_p}\overline{S}\ar[d]^-{c\lambda_0\otimes\id}\ar[rr]^-{\rho_L}&&
	    X_L\times_S\overline{S}\ar[d]^-{p_L\times\id}_-{\cong}\\
	    \X^\vee\otimes_{\overline{\F}_p}\overline{S}&&
	    (X_{L^\vee})^\vee\times_S\overline{S}\ar[ll]_-{\rho^\vee_{L^\vee}}\lefteqn{.}
	    }
	    \]
 \item[(d)] For each $L\in\mathscr{L}$, we have the following equality
	    of polynomial functions on $a\in \mathcal{O}_B$:
	    \[
	    \det\nolimits_{\mathcal{O}_S}(a;\Lie X_L)=\det\nolimits_K(a;V_0).
	    \]
	    This equation is called the determinant condition. 
	    For a precise formulation, see \cite[3.23 a)]{MR1393439}.
\end{enumerate}
 The functor $\M$ is represented by a formal scheme (denoted by the same symbol $\M$) which is
 locally formally of finite type over $\Spf \mathcal{O}_{\breve{E}}$ (see \cite[Theorem 3.25]{MR1393439}).
 We write $\overline{\M}$ for the underlying reduced scheme $\M^{\mathrm{red}}$ of $\M$. It is a scheme locally of finite type
 over $\overline{\F}_p$, which is often not quasi-compact. Each irreducible component of $\overline{\M}$ is known to be projective
 over $\overline{\F}_p$ (see \cite[Proposition 2.32]{MR1393439}).

 The group $J=\mathbf{J}(\Q_p)$ acts naturally on the functor $\M$ on the left; the element $h\in J$, regarded as a self-quasi-isogeny
 on $\X$, carries $\{(X_L,\iota_L,\rho_L)\}_{L\in\mathscr{L}}$ to 
 $\{(X_L,\iota_L,\rho_L\circ h^{-1})\}_{L\in\mathscr{L}}$. Hence $J$ also acts on the formal scheme $\M$ and the scheme $\overline{\M}$.
\end{defn}

\begin{rem}
 In the definition of $\M$, we fixed a triple $(\X,\lambda_0,\iota_0)$. 
 However, we can check that the formal scheme $\M$ with the action of $J$ is essentially independent of this choice.
 See the remark after \cite[Definition 3.21]{MR1393439}.
\end{rem}

For a scheme $X$ locally of finite type over $\overline{\F}_p$, we write $\Irr(X)$ for the set of
irreducible components.
We investigate the action of $J$ on $\Irr(\overline{\M})$ under the following condition:

\begin{defn}\label{defn:comes-from-AV}
 We say that $b$ comes from an abelian variety if there exist
 \begin{itemize}
  \item a $\Q$-subalgebra $\widetilde{B}$ which is stable under $*$ and satisfies $\widetilde{B}\otimes_\Q\Q_p=B$,
  \item an abelian variety $A_0$ over $\overline{\F}_p$,
  \item a polarization $\lambda_0\colon A_0\longrightarrow A_0^\vee$,
  \item and a homomorphism $\iota_0\colon \widetilde{B}\longrightarrow \End(A_0)\otimes_\Z\Q$
	satisfying $\lambda_0\circ\iota_0(a^*)=\iota_0(a)^\vee\circ\lambda_0$ for every $a\in \widetilde{B}$
 \end{itemize} 
 such that the polarized $B$-isocrystal $\mathbb{D}(A_0[p^\infty])_{\Q}$ is isomorphic to $N_b$.

 Note that in this case $\widetilde{B}$ is a finite-dimensional semisimple algebra over $\Q$,
 $*$ is a positive involution on $\widetilde{B}$,
 and $\iota_0$ is an injective homomorphism (recall that we assume $V$ to be a faithful $B$-module).
\end{defn}

\begin{rem}\label{rem:Manin}
 The problem of finding $(A_0,\lambda_0,\iota_0)$ as in Definition \ref{defn:comes-from-AV}
 for a fixed $\widetilde{B}$
 can be regarded as a generalization of Manin's problem.
 In fact, the case where $(B,*,\widetilde{B})=(\Q_p,\id,\Q)$ is exactly the original problem
 (\cite[Chapter IV, \S 5, Conjecture 1, Conjecture 2]{MR0157972}), which has been solved in
 \cite{MR3077121} and \cite{MR1792294}. In the next subsection, by using a similar method as in \cite{MR3077121},
 we prove that any $b$ comes from an abelian variety if $\widetilde{B}$ is a quaternion division algebra over $\Q$
 which ramifies exactly at $p$ and another finite place.

 This problem is studied by many people in the context of non-emptiness of
 Newton strata of Shimura varieties;
 see \cite{MR2789744}, \cite{MR3072810}, \cite{MR2983012}, \cite{Kret-PEL} for example. 
 Let us assume that there exists
 \begin{itemize}
  \item a semisimple $\Q$-algebra $\widetilde{B}$,
  \item an order $\mathcal{O}_{\widetilde{B}}$ of $\widetilde{B}$,
  \item a positive involution $\widetilde{*}$ on $\widetilde{B}$ under which $\mathcal{O}_{\widetilde{B}}$ is stable,
  \item a finite faithful $\widetilde{B}$-module $\widetilde{V}$,
  \item and a non-degenerate alternating bilinear pairing $\langle\ ,\ \rangle^\sim\colon \widetilde{V}\times \widetilde{V}\longrightarrow \Q$
	satisfying $\langle av,w\rangle^\sim=\langle v,a^{\widetilde{*}}w\rangle^\sim$ for every $a\in \widetilde{B}$
	and $v,w\in \widetilde{V}$
 \end{itemize} 
 such that the localization of $(\widetilde{B},\widetilde{*},\widetilde{V},\langle\ ,\ \rangle^\sim)$ at $p$
 is isomorphic to $(B,*,V,\langle\ ,\ \rangle)$. 
 We write $\widetilde{\mathbf{G}}$ for the algebraic group over $\Q$
 consisting of $\widetilde{B}$-linear automorphisms of $\widetilde{V}$ which preserve $\langle\ ,\ \rangle^\sim$
 up to a scalar multiple. We fix an isomorphism $\C\cong \overline{K}_0$ and assume further that
 \begin{itemize}
  \item there exists a Shimura datum $(\widetilde{\mathbf{G}},h)$ (in particular $\widetilde{\mathbf{G}}$ is
	connected) such that the associated Hodge cocharacter $\mu_h\colon \mathbb{G}_{m,\C}\longrightarrow \widetilde{\mathbf{G}}_\C$ is identified with $\mu\colon \mathbb{G}_{m,\overline{K}_0}\longrightarrow \mathbf{G}_{\overline{K}_0}\cong\widetilde{\mathbf{G}}_{\overline{K}_0}$ under the fixed isomorphism $\C\cong \overline{K}_0$.
 \end{itemize}
 Finally, we assume that $(\mu,b)$ is admissible in the sense of \cite[Definition 1.18]{MR1393439}.
 This condition is usually included in the definition of Rapoport-Zink data
 (see \cite[Definition 3.18]{MR1393439}).
In this situation, every $b$ comes from an abelian variety if the datum $(B,\mathcal{O}_B,*,V,\langle\ ,\ \rangle)$ is unramified, namely, $B$ is a product of matrix algebras over unramified extensions of $\Q_p$
and there exists a $\Z_p$-lattice of $V$ which is self-dual for $\langle\ ,\ \rangle$
and preserved by $\mathcal{O}_B$. See \cite[Theorem 10.1]{MR3072810}.
\end{rem}

The main theorem of this paper is the following:

\begin{thm}\label{thm:main}
 Assume that $b$ comes from an abelian variety.
 Then, $\Irr(\overline{\M})$ has finitely many $J$-orbits.
\end{thm}

\begin{rem}\label{rem:known-cases}
 \begin{enumerate}
  \item If the datum $(B,\mathcal{O}_B,*,V,\langle\ ,\ \rangle)$ is unramified,
	the corresponding result is obtained in \cite[Th\'eor\`eme 2.4.13]{MR2074714}.
	Together with Remark \ref{rem:Manin}, the theorem above gives an alternative proof of Fargues' result.
  \item In some concrete cases, there are more precise results on irreducible components of $\overline{\M}$. 
	See \cite{MR2369090}, \cite{MR2466182}, \cite{MR2800696}, \cite{MR3175176} and \cite{MR3272278}
	for instance.
 \end{enumerate}
\end{rem}

\begin{rem}
 In this paper, we only work on Rapoport-Zink spaces of PEL type. 
 However, it will be possible to apply the same technique to Rapoport-Zink spaces of EL type.
\end{rem}

\subsection{Example: quaternion unitary case}\label{subsec:GU(n,D)}
In this subsection, we will give an example of $(B,*,V,\langle\ ,\ \rangle)$ such that
every $b$ as in the previous subsection comes from an abelian variety.
Let $D$ be a quaternion division algebra over $\Q_p$.
Recall that $D$ can be written as $\Q_{p^2}[\Pi]$,
where $\Q_{p^2}$ is the unramified quadratic extension of $\Q_p$,
and $\Pi$ is an element satisfying $\Pi^2=p$ and $\Pi a=\sigma(a)\Pi$ for every $a\in \Q_{p^2}$ and
the unique non-trivial element $\sigma$ of $\Gal(\Q_{p^2}/\Q_p)$.
Let $*$ be an involution of $D$ defined by $\Pi^*=\Pi$ and $a^*=\sigma(a)$ for every $a\in \Q_{p^2}$.
Note that we have $d^*=\varepsilon (\Trd(d)-d) \varepsilon^{-1}$
for any $\varepsilon\in \Q_{p^2}^\times$ with $\sigma(\varepsilon)=-\varepsilon$,
where $\Trd$ denotes the reduced trace.
We fix an integer $n\ge 1$ and a non-degenerate alternating bilinear pairing
$\langle\ ,\ \rangle\colon D^n\times D^n\longrightarrow \Q_p$
satisfying $\langle dx,y\rangle=\langle x,d^*y\rangle$ for every $d\in D$ and $x,y\in D^n$.
An example of such a pairing $\langle\ ,\ \rangle$ is given by
$\langle (d_i),(d'_i)\rangle=\sum_{i=1}^n\Trd(\varepsilon d_i^*d'_i)$,
where $\varepsilon\in \Q_{p^2}^\times$ is an element satisfying $\sigma(\varepsilon)=-\varepsilon$.
For the quadruple $(D,*,D^n,\langle\ ,\ \rangle)$, the algebraic group $\mathbf{G}$ introduced in the previous subsection
is called the quaternion unitary similitude group.

In this subsection, we prove the following proposition.

\begin{prop}\label{prop:GU(n,D)-comes-from-AV}
 Let $b$ be an element of $\mathbf{G}(K_0)$ satisfying the following conditions:
 \begin{itemize}
  \item the isocrystal $(N_b,\Phi_b)$ has slopes in the interval $[0,1]$,
  \item and $\simil(b)=p$.
 \end{itemize}
 Then, $b$ comes from an abelian variety.
\end{prop}

We begin with describing polarized $D$-isocrystals.

\begin{lem}\label{lem:D-decomposition}
 \begin{enumerate}
  \item For a $\Q_{p^2}$-algebra $R$, the following two categories are equivalent:
	\begin{itemize}
	 \item The category of pairs $(N,\langle\ ,\ \rangle)$ consisting of an $R$-module $N$ endowed with
	       an action of $D$ and a non-degenerate alternating bilinear pairing
	       $\langle\ ,\ \rangle\colon N\times N\longrightarrow R$ satisfying
	       $\langle dx,y\rangle=\langle x,d^*y\rangle$ for every $d\in D$ and $x,y\in N$.
	       The morphisms are $R$-homomorphisms preserving the $D$-actions and
	       the $R^\times$-homothety classes of the pairings.
	 \item The category of pairs $(N_0,\langle\ ,\ \rangle_0)$ consisting of an $R$-module $N_0$
	       and a non-degenerate alternating bilinear pairing
	       $\langle\ ,\ \rangle_0\colon N_0\times N_0\longrightarrow R$.
	       The morphisms are $R$-homomorphisms preserving the $R^\times$-homothety classes of the pairings.
	\end{itemize}
  \item The following two categories are equivalent:
	\begin{itemize}
	 \item The category of polarized $D$-isocrystals.
	 \item The category of polarized $\Q_p$-isocrystals, where we consider the identity as an involution on $\Q_p$.
	\end{itemize}
 \end{enumerate}
\end{lem}

\begin{prf}
 This lemma is essentially obtained in \cite[1.42]{MR1393439}, but we include the proof here
 in order to fix notation.

 We prove i).
 Let $(N,\langle\ ,\ \rangle)$ be an object of the first category.
 For $i=0,1$, we denote by $N_i$ the $R$-submodule of $N$ on which $\Q_{p^2}\subset D$ acts by
 $\Q_{p^2}\yrightarrow{\sigma^i}\Q_{p^2}\longrightarrow R$. Then we have $N=N_0\oplus N_1$.
 Note that the action of $\Pi\in D$ induces an isomorphism between $N_0$ and $N_1$, hence $N\cong N_0\oplus N_0$.
 We can easily check that $N_0$ and $N_1$ are totally isotropic with respect to $\langle\ ,\ \rangle$.

 Let $\langle\ ,\ \rangle_0\colon N_0\times N_0\longrightarrow R$ be the pairing given by
 $\langle x,y\rangle_0=\langle x,\Pi y\rangle$.
 It is non-degenerate and alternating.
 Therefore we obtain the functor $(N,\langle\ ,\ \rangle)\longmapsto (N_0,\langle\ ,\ \rangle_0)$ from the first category
 to the second. Note also that we have
 $\langle x+\Pi x',y+\Pi y'\rangle=\langle x,y'\rangle_0+\langle x',y\rangle_0$ for $x,x',y,y'\in N_0$.

 We shall construct a quasi-inverse functor. Let $(N_0,\langle\ ,\ \rangle_0)$ be an object of the second category.
 We define an action of $D$ on $N_0\oplus N_0$ by $\Pi\cdot (x,y)=(py,x)$ and $a\cdot (x,y)=(ax,\sigma(a)y)$ for $a\in \Q_{p^2}$.
 Let $\langle\ ,\ \rangle\colon (N_0\oplus N_0)\times (N_0\oplus N_0)\longrightarrow R$ be a pairing
 defined by $\langle (x_0,x_1),(y_0,y_1)\rangle=\langle x_0,y_1\rangle_0+\langle x_1,y_0\rangle_0$
 for $x_0,x_1,y_0,y_1\in N_0$. 
 It is easy to observe that the functor $(N_0,\langle\ ,\ \rangle_0)\longmapsto (N_0\oplus N_0,\langle\ ,\ \rangle)$ is a quasi-inverse of $(N,\langle\ ,\ \rangle)\longmapsto (N_0,\langle\ ,\ \rangle_0)$. This concludes i).

 Next we consider ii). Let $(N,\Phi,\langle\ ,\ \rangle,c)$ be a polarized $D$-isocrystal.
 Then $\Phi$ maps $N_0$ to $N_1$, hence $\Pi^{-1}\Phi\colon N_0\longrightarrow N_0$ makes $N_0$ an isocrystal.
 Since
 \[
 \langle \Pi^{-1}\Phi x,\Pi^{-1}\Phi y\rangle_0=\langle \Pi^{-1}\Phi x,\Phi y\rangle
 =c\sigma(\langle \Pi^{-1}x,y\rangle)
 =c\sigma(\langle x,\Pi^{-1}y\rangle)
 =p^{-1}c\sigma(\langle x,y\rangle_0)
 \]
 for $x,y\in N_0$, $(N_0,\Pi^{-1}\Phi,\langle\ ,\ \rangle_0,p^{-1}c)$ is a polarized $\Q_p$-isocrystal.

 Conversely, let $(N_0,\Phi_0,\langle\ ,\ \rangle_0,c_0)$ be a polarized $\Q_p$-isocrystal.
 Let $\langle\ ,\ \rangle\colon (N_0\oplus N_0)\times (N_0\oplus N_0)\longrightarrow K_0$
 be the pairing constructed in i).
 We define $\Phi\colon N_0\oplus N_0\longrightarrow N_0\oplus N_0$ by $(x,y)\longmapsto (p\Phi_0 y,\Phi_0 x)$.
 Then $(N_0\oplus N_0,\Phi,\langle\ ,\ \rangle,pc_0)$ is a polarized $D$-isocrystal,
 and $(N_0,\Phi_0,\langle\ ,\ \rangle_0,c_0)\longmapsto (N_0\oplus N_0,\Phi,\langle\ ,\ \rangle,pc_0)$ gives a quasi-inverse
 of the functor $(N,\Phi,\langle\ ,\ \rangle,c)\longmapsto (N_0,\Pi^{-1}\Phi,\langle\ ,\ \rangle_0,p^{-1}c)$.
 This concludes ii).
\end{prf}

\begin{rem}\label{rem:D-decomp}
 \begin{enumerate}
  \item Set $R=K_0$ in Lemma \ref{lem:D-decomposition} i).
	For an integer $d\ge 1$, all objects $(N_0,\langle\ ,\ \rangle_0)$ in the second category of Lemma \ref{lem:D-decomposition} i)
	with $\dim_{K_0}N_0=2d$ are isomorphic.
	Therefore, all objects $(N,\langle\ ,\ \rangle)$ in the first category of Lemma \ref{lem:D-decomposition} i)
	with $\dim_{K_0}N=4d$ are isomorphic.
	This implies that each polarized $D$-isocrystal $(N,\Phi,\langle\ ,\ \rangle,c)$ with $\dim_{K_0}N=4n$ is
	isomorphic to $(N_b,\Phi_b,\langle\ ,\ \rangle,\simil(b))$ for some $b\in \mathbf{G}(K_0)$.
  \item Lemma \ref{lem:D-decomposition} i) gives an inner twist
	$\mathbf{G}\otimes_{\Q_p}\Q_{p^2}\yrightarrow{\cong}\GSp_{2n}\otimes_{\Q_p}\Q_{p^2}$.
 \end{enumerate}
\end{rem}

\begin{defn}
 The Newton polygon of a polarized $D$-isocrystal $(N,\Phi,\langle\ ,\ \rangle,c)$
 is defined to be the Newton polygon of the isocrystal $(N_0,\Pi^{-1}\Phi)$.
\end{defn}

\begin{rem}\label{rem:GL-Newton}
 Let $(N,\Phi,\langle\ ,\ \rangle,c)$ be a polarized $D$-isocrystal, and $\sum_{i=0}^r m_i[\lambda_i]$ its Newton polygon.
 Here $\lambda_0<\cdots<\lambda_r$ denote the slopes, and $m_i$ denotes the multiplicity of the slope $\lambda_i$ part.
 Then we have $\lambda_{r-i}=\ord_p(c)-1-\lambda_i$ and $m_{r-i}=m_i$ for every $0\le i\le r$,
 where $\ord_p$ denotes the $p$-adic order.
 Further, the Newton polygon of the isocrystal $(N,\Phi)$ is given by $\sum_{i=0}^r m'_i[\lambda_i+1/2]$,
 where $m'_i=2m_i$ if $\lambda_i$ and $\lambda_i+1/2$ have the same denominator, and $m'_i=m_i$ otherwise.
 In particular, the slopes of $(N,\Phi)$ lie in the interval $[0,1]$ if and only if $-1/2\le \lambda_i\le 1/2$ for every $i$.
\end{rem}

\begin{lem}\label{lem:D-isoc-Newton-inj}
 Two polarized $D$-isocrystals are isomorphic if and only if the corresponding Newton polygons are the same.
\end{lem}

\begin{prf}
 We write $\mathbf{B}(\mathbf{G})$ for the set of $\sigma$-conjugacy classes of $\mathbf{G}(K_0)$.
 As in \cite[\S 4.2]{MR1485921}, the inner twist $\mathbf{G}\otimes_{\Q_p}\Q_{p^2}\yrightarrow{\cong}\GSp_{2n}\otimes_{\Q_p}\Q_{p^2}$
 in Remark \ref{rem:D-decomp} ii) determines the Newton map $\nu_{\mathbf{G}}\colon \mathbf{B}(\mathbf{G})\longrightarrow \R^n$.
 As the derived group of $\mathbf{G}$ is simply connected,
 $\nu_{\mathbf{G}}$ is known to be injective (see \cite[\S 4.13]{MR1485921}).

 By Remark \ref{rem:D-decomp} i), the map sending $b$ to $(N_b,\Phi_b,\langle\ ,\ \rangle,\simil(b))$
 gives a bijection between $\mathbf{B}(\mathbf{G})$ and the set of isomorphism classes of polarized $D$-isocrystals
 whose dimensions as $K_0$-vector spaces are $4n$.
 By definition, for $b\in \mathbf{B}(\mathbf{G})$, the Newton polygon of the polarized
 $D$-isocrystal $(N_b,\Phi_b,\langle\ ,\ \rangle,\simil(b))$ has essentially the same information as $\nu_{\mathbf{G}}(b)$.
 Hence the injectivity of $\nu_{\mathbf{G}}$ provides the desired claim.
\end{prf}

\begin{defn}\label{defn:simple-D-isoc}
 For integers $m\ge 0$, $k>0$ with $(m,k)=1$,
 let $N_{m,k}$ denote a $p$-polarized $D$-isocrystal whose Newton polygon equals $[-m/k]+[m/k]$.
 Such an $N_{m,k}$ exists by Lemma \ref{lem:D-decomposition} ii), and it is unique up to isomorphism
 by Lemma \ref{lem:D-isoc-Newton-inj}.
 We call $N_{m,k}$ the simple polarized $D$-isocrystal of type $(m,k)$.
\end{defn}

\begin{lem}\label{lem:decomp-into-simples}
 Every $p$-polarized $D$-isocrystal is isomorphic to a direct sum of simple polarized $D$-isocrystals of various types.
\end{lem}

\begin{prf}
 Let $N$ be a $p$-polarized $D$-isocrystal.
 By Remark \ref{rem:GL-Newton}, it is easy to find simple polarized $D$-isocrystals
 $N_{m_1,k_1},\ldots,N_{m_s,k_s}$ such that the Newton polygon of
 $N_{m_1,k_1}\oplus \cdots\oplus N_{m_s,k_s}$ is equal to that of $N$.
 By Lemma \ref{lem:D-isoc-Newton-inj}, we conclude that $N$ is isomorphic to $N_{m_1,k_1}\oplus \cdots\oplus N_{m_s,k_s}$.
\end{prf}

We fix a prime number $q$ different from $p$.
Let $\widetilde{D}$ be the quaternion algebra over $\Q$ ramified exactly at $p$ and $q$.

\begin{lem}
 There exists an embedding $\widetilde{D}\hooklongrightarrow D$ such that
 $\widetilde{D}$ is stable under $*$ and the restriction of $*$ to $\widetilde{D}$ is a positive involution.
\end{lem}

\begin{prf}
 We can find an element $\varepsilon\in \widetilde{D}^\times$ such that $\varepsilon^2\in \Q^\times$ and
 $\Q(\varepsilon)$ is an imaginary quadratic field which is non-split and unramified at $p$.
 By the Skolem-Noether theorem, there exists an isomorphism $\widetilde{D}\otimes_\Q\Q_p\yrightarrow{\cong}D$
 which maps $\varepsilon$ to an element of $\Q_{p^2}$.
 We regard $\widetilde{D}$ as a $\Q$-subalgebra of $D$ by this isomorphism.
 As $\varepsilon\in \Q^\times_{p^2}$ satisfies $\sigma(\varepsilon)=-\varepsilon$, 
 we have $d^*=\varepsilon (\Trd(d)-d)\varepsilon^{-1}$.
 This implies that $\widetilde{D}$ is stable under $*$.
 Since $\varepsilon^2<0$, the restriction of $*$ to $\widetilde{D}$ is a positive involution
 (see \cite[\S 21, Theorem 2]{MR0282985}).
\end{prf}

We fix an embedding $\widetilde{D}\hooklongrightarrow D$ as in the previous lemma.
Thanks to Lemma \ref{lem:decomp-into-simples} and Remark \ref{rem:GL-Newton},
Proposition \ref{prop:GU(n,D)-comes-from-AV} is reduced to the following proposition:

\begin{prop}\label{prop:simple-case}
 For integers $m\ge 0$, $n>0$ such that $(m,n)=1$ and $0\le 2m\le n$, there exists
 \begin{itemize}
  \item a $2n$-dimensional abelian variety $A_0$ over $\overline{\F}_p$,
  \item a polarization $\lambda_0\colon A_0\longrightarrow A_0^\vee$,
  \item and a homomorphism $\iota_0\colon \widetilde{D}\longrightarrow \End(A_0)\otimes_\Z\Q$
	satisfying $\lambda_0\circ\iota_0(a^*)=\iota_0(a)^\vee\circ\lambda_0$ for every $a\in \widetilde{D}$
 \end{itemize}
 such that the $p$-polarized $D$-isocrystal $\mathbb{D}(A_0[p^\infty])_\Q$ is isomorphic to $N_{m,n}$.
\end{prop}

\begin{prf}
 We follow the method in \cite[\S 1, Un exemple sp\'ecial]{MR3077121}.

 First we assume that $n$ is odd and greater than $1$. Note that in this case $m\ge 1$.
 We choose an integer $u$ such that $p\nmid u$ and $x^2+p^{n-2m}x+p^{2n}u$ is irreducible
 in $\F_q$.
 Let $\pi\in\overline{\Q}$ be a root of the equation $x^2+p^{n-2m}x+p^{2n}u=0$ 
 and $A_0$ the abelian variety over $\F_{p^{2n}}$
 corresponding to $\pi$ under the bijection in \cite[Th\'eor\`eme 1 (i)]{MR3077121}. 
 We put $F=\Q(\pi)$. It is an imaginary quadratic field in which $p$ splits into two places $v$, $\overline{v}$
 such that $v(\pi)=n-2m$ and $\overline{v}(\pi)=n+2m$.
 By the choice of $u$, $F$ does not split at $q$.
 By \cite[Th\'eor\`eme 1 (ii)]{MR3077121}, the endomorphism algebra $\End(A_0)\otimes_\Z\Q$ is a central simple
 algebra over $F$ whose invariants are given by
 \[
  \inv_w (\End(A_0)\otimes_\Z\Q)=\begin{cases}
				  0& w\nmid p,\\
				  \frac{n-2m}{2n}& w=v,\\
				  \frac{n+2m}{2n}& w=\overline{v},
				 \end{cases}
 \]
 where $w$ is a place of $F$. In particular we have $\dim_F(\End(A_0)\otimes_\Z\Q)=(2n)^2$
 and $\dim A_0=2n$.
 We can observe that $\inv_w (\End(A_0)\otimes_\Z\Q)-\inv_w (\widetilde{D}\otimes_\Q F)$ is killed by $n$
 for every place $w$ of $F$; note that if $w\mid q$, 
 then $\widetilde{D}\otimes_\Q F$ splits at $w$ since $F$ does not split at $q$.
 Therefore, there exists an embedding 
 $\widetilde{D}\otimes_\Q F\hooklongrightarrow \End(A_0)\otimes_\Z\Q$.
 Let $\iota_0$ be the composite of 
 $\widetilde{D}\hooklongrightarrow\widetilde{D}\otimes_\Q F\hooklongrightarrow \End(A_0)\otimes_\Z\Q$.
 By \cite[Lemma 9.2]{MR1124982}, there exists a polarization $\lambda_0\colon A_0\longrightarrow A_0^\vee$
 satisfying $\lambda_0\circ\iota_0(a^*)=\iota_0(a)^\vee\circ\lambda_0$ for every $a\in \widetilde{D}$.
 By construction, the Newton polygons of
 $\mathbb{D}((A_0\otimes_{\F_{p^{2n}}}\overline{\F}_p)[p^\infty])_\Q$ and $N_{m,n}$
 regarded as isocrystals are the same.
 By Remark \ref{rem:GL-Newton}, the Newton polygons of
 $\mathbb{D}((A_0\otimes_{\F_{p^{2n}}}\overline{\F}_p)[p^\infty])_\Q$ and $N_{m,n}$
 regarded as polarized $D$-isocrystals are also the same.
 Hence, by Lemma \ref{lem:D-isoc-Newton-inj} we conclude that the $p$-polarized $D$-isocrystal 
 $\mathbb{D}((A_0\otimes_{\F_{p^{2n}}}\overline{\F}_p)[p^\infty])_\Q$ is isomorphic to $N_{m,n}$.

 Next we assume that $n$ is even. 
 Let $\pi\in\overline{\Q}$ be a root of the equation $x^2+p^{\frac{n}{2}-m}x+p^n=0$ 
 and $A_1$ the abelian variety over $\F_{p^n}$
 corresponding to $\pi$. We put $A_0=A_1\times A_1$. Then, by the same argument as above,
 we have $\dim A_0=2n$ and we can find a homomorphism
 $\iota_0\colon \widetilde{D}\longrightarrow \End(A_0)\otimes_\Z\Q=M_2(\End(A_1)\otimes_\Z\Q)$
 and a polarization $\lambda_0\colon A_0\longrightarrow A_0^\vee$.
 We can observe that $\mathbb{D}((A_0\otimes_{\F_{p^n}}\overline{\F}_p)[p^\infty])_\Q\cong N_{m,n}$
 by comparing their Newton polygons.
 
 Finally we consider the case $n=1$. Note that $m=0$.
 Let $E_0$ be a supersingular elliptic curve over $\overline{\F}_p$ and put $A_0=E_0\times E_0$.
 Then $\End(E_0)\otimes_\Z\Q$ is a quaternion algebra over $\Q$ which ramifies exactly at $\infty$ and $p$,
 and $\End(A_0)\otimes_\Z\Q$ equals $M_2(\End(E_0)\otimes_\Z\Q)$.
 We can immediately check that there exists a homomorphism
 $\iota_0\colon \widetilde{D}\hooklongrightarrow \End(A_0)\otimes_\Z\Q$.
 By using it, we can repeat the argument above.
\end{prf}

Now the proof of Proposition \ref{prop:GU(n,D)-comes-from-AV} is complete.

\section{Proof of the main theorem}\label{sec:proof-main-thm}
\subsection{First reduction}
In the first step of our proof of Theorem \ref{thm:main}, we will replace $\M$ with
a simpler moduli space $\mathcal{N}$.
We fix a lattice $L\in \mathscr{L}$ such that $L\subset L^\vee$.

\begin{defn}\label{defn:N}
 Let $\mathcal{N}\colon \Nilp_{W(\overline{\F}_p)}\longrightarrow \mathbf{Set}$ be
 the functor that associates $S$ with
 the set of isomorphism classes of $(X,\iota,\rho)$ where
 \begin{itemize}
  \item $X$ is a $p$-divisible group over $S$,
  \item $\iota\colon \mathcal{O}_B\longrightarrow \End(X)$ is a homomorphism,
  \item and $\rho\colon \X\otimes_{\overline{\F}_p}\overline{S}\longrightarrow X\times_S\overline{S}$ is a quasi-isogeny
	compatible with $\mathcal{O}_B$-actions,
 \end{itemize}
 such that the following condition is satisfied.
 \begin{itemize}
  \item Locally on $S$, there exist a constant $c\in \Q_p^\times$ and
	an isogeny $\lambda\colon X\longrightarrow X^\vee$ with $\height\lambda=\log_p\#(L^\vee/L)$
	such that the following diagram is commutative:
	\[
	\xymatrix{%
	\X\otimes_{\overline{\F}_p}\overline{S}\ar[d]^-{c\lambda_0\otimes\id}\ar[rr]^-{\rho}&&
	X\times_S\overline{S}\ar[d]^-{\lambda\times\id}\\
	\X^\vee\otimes_{\overline{\F}_p}\overline{S}&&
	X^\vee\times_S\overline{S}\ar[ll]_-{\rho^\vee}\lefteqn{.}
	}
	\]
 \end{itemize}
 By a similar method as in the proof of \cite[Theorem 3.25]{MR1393439}, we can prove that $\mathcal{N}$ is represented by
 a formal scheme (denoted by the same symbol $\mathcal{N}$) which is locally formally of finite type over $\Spf W(\overline{\F}_p)$.
 We write $\overline{\mathcal{N}}$ for the underlying reduced scheme $\mathcal{N}^{\mathrm{red}}$
 of $\mathcal{N}$.
 
 The group $J$ acts naturally on $\mathcal{N}$ on the left;
 an element $h\in J$ carries $(X,\iota,\rho)$ to $(X,\iota,\rho\circ h^{-1})$.
\end{defn}

\begin{lem}\label{lem:M-N-proper}
 There exists a proper morphism of formal schemes
 $\M\longrightarrow \mathcal{N}$ compatible with $J$-actions.
\end{lem}

\begin{prf}
 It suffices to construct a proper morphism $\M\longrightarrow \mathcal{N}\otimes_{W(\overline{\F}_p)}\mathcal{O}_{\breve{E}}$.
 Let $S$ be an object of $\Nilp_{\mathcal{O}_{\breve{E}}}$ and $\{(X_{L'},\iota_{L'},\rho_{L'})\}_{L'\in\mathcal{L}}$ an element of $\M(S)$. We will show that $(X_L,\iota_L,\rho_L)$ gives an element of
$\mathcal{N}(S)$. We need to verify the condition on quasi-polarizations.
 By the condition (c) in Definition \ref{defn:RZ-space}, locally on $S$ there exist $c\in \Q_p^\times$
 and an isomorphism $p_L\colon X_L\longrightarrow (X_{L^\vee})^\vee$.
 Let $\lambda$ be the composite $X_L\yrightarrow[\cong]{p_L}(X_{L^\vee})^\vee\yrightarrow{(\widetilde{\rho}_{L^\vee,L})^\vee}X_L^\vee$. Then, we have 
 \begin{gather*}
  \height\lambda=\height\widetilde{\rho}_{L^\vee,L}=\log_p\#(L^\vee/L),\\
  c\lambda_0\otimes\id=\rho^\vee_{L^\vee}\circ (p_L\times\id)\circ \rho_L=\rho_L^\vee\circ (\lambda\times\id)\circ \rho_L,
 \end{gather*}
 as desired. Hence we obtain a morphism $\M\longrightarrow \mathcal{N}\otimes_{W(\overline{\F}_p)}\mathcal{O}_{\breve{E}}$ of formal schemes, which we denote by $\psi$. Clearly $\psi$ is $J$-equivariant.

 Next we prove that $\psi$ is proper. We denote the universal object over $\mathcal{N}\otimes_{W(\overline{\F}_p)}\mathcal{O}_{\breve{E}}$ by $(\widetilde{X},\widetilde{\iota},\widetilde{\rho})$.
 We want to classify $\{(X_{L'},\iota_{L'},\rho_{L'})\}_{L'\in\mathscr{L}}$ as 
 in Definition \ref{defn:RZ-space} such that 
 $(X_L,\iota_L,\rho_L)=(\widetilde{X},\widetilde{\iota},\widetilde{\rho})$.
 By the condition (b) in Definition \ref{defn:RZ-space}, such an object is determined by
 $\{(X_{L'},\iota_{L'},\rho_{L'})\}_{L'\in\mathscr{L},L\subset L'\subsetneq p^{-1}L}$
 up to isomorphism. 
 Furthermore, for each $L'\in\mathscr{L}$ with $L\subset L'\subsetneq p^{-1}L$,
 $(X_{L'},\iota_{L'},\rho_{L'})$ is determined by $\Ker\widetilde{\rho}_{L',L}$, which is
 a finite flat subgroup scheme of $X_L[p]$ with degree $\#(L'/L)$.
 Clearly such a finite flat subgroup scheme is classified by a formal scheme $\M_{L'}$
 which is proper over $\mathcal{N}\otimes_{W(\overline{\F}_p)}\mathcal{O}_{\breve{E}}$.
 Let $\M'$ be the fiber product over $\mathcal{N}\otimes_{W(\overline{\F}_p)}\mathcal{O}_{\breve{E}}$
 of $\M_{L'}$ for $L'\in\mathscr{L}$ with $L\subset L'\subsetneq p^{-1}L$.
 By \cite[Proposition 2.9]{MR1393439}, it is easy to observe that the natural morphism
 $\M\longrightarrow \M'$ is a closed immersion
 (strictly speaking, we use the fact that \cite[Proposition 2.9]{MR1393439} is valid for
 $\alpha\in\Hom(X,Y)\otimes_{\Z_p}\Q_p$ which is not necessary a quasi-isogeny).
 This concludes the proof.
\end{prf}

\begin{cor}
 If the action of $J$ on $\Irr(\overline{\mathcal{N}})$ has finite orbits,
 then so does the action of $J$ on $\Irr(\overline{\M})$.
\end{cor}

\begin{prf}
 Take finitely many $\alpha_1,\ldots,\alpha_m\in \Irr(\overline{\mathcal{N}})$
 such that $\Irr(\overline{\mathcal{N}})=\bigcup_{i=1}^m J\alpha_i$.
 By Lemma \ref{lem:M-N-proper}, for each $i$, only finitely many components 
 $\beta_1,\ldots,\beta_{k_i}\in \Irr(\overline{\mathcal{M}})$ are mapped into $\alpha_i$ by the morphism 
 $\overline{\M}\longrightarrow \overline{\mathcal{N}}$. It is easy to observe that
 $\Irr(\overline{\M})=\bigcup_{i=1}^m\bigcup_{j=1}^{k_i}J\beta_j$.
\end{prf}

By this corollary, we may consider $\overline{\mathcal{N}}$ in place of $\overline{\M}$.

\begin{rem}
 The determinant condition (see Definition \ref{defn:RZ-space} (d)) plays an important role 
 when we relate the moduli space $\mathcal{M}$ with Shimura varieties.
 If we remove it from the conditions defining the functor $\mathcal{M}$, we obtain a bigger moduli space 
 $\mathcal{M}'$ containing $\mathcal{M}$ as a closed formal subscheme.
 It is (at least a priori) nothing to do with any Shimura variety, but still useful for our purpose
 since $\Irr(\overline{\mathcal{M}'})$ gives an ``upper bound'' of $\Irr(\overline{\mathcal{M}})$
 in some sense.
 Our new moduli space $\mathcal{N}$ is a slight modification of $\mathcal{M}'$.
\end{rem}

\subsection{A moduli space of PEL type}
So far, we considered an arbitrary triple $(\X,\lambda_0,\iota_0)$ as in Section \ref{sec:RZ-space}. Now we take $(A_0,\lambda_0,\iota_0)$
as in Definition \ref{defn:comes-from-AV} and consider the triple $(\X,\lambda_0,\iota_0)$
attached to it; namely, we put $\X=A_0[p^\infty]$ and denote the induced quasi-polarization $\X\longrightarrow \X^\vee$ 
and homomorphism $B\longrightarrow \End(\X)\otimes_{\Z_p}\Q_p$ by the same symbols $\lambda_0$ and $\iota_0$.

\begin{lem}\label{lem:replace-AV}
 Assume that $\mathcal{N}(\overline{\F}_p)\neq \varnothing$. 
 Then we can replace $(A_0,\lambda_0,\iota_0)$ so that the following conditions are satisfied:
 \begin{itemize}
  \item $\X=A_0[p^\infty]$ is completely slope divisible (\cf \cite[Definition 10]{MR1844205}).
  \item $\ord_p\deg\lambda_0=\log_p\#(L^\vee/L)$.
  \item there exists an order $\widetilde{\mathcal{O}}$ of $\widetilde{B}$ which is
	contained in $\iota_0^{-1}(\End(A_0))$, stable under $*$
	and satisfies $\widetilde{\mathcal{O}}\otimes_\Z\Z_p=\mathcal{O}_B$.
 \end{itemize}
\end{lem}

\begin{prf}
 First we prove that $\mathcal{N}(\overline{\F}_p)$ contains an element $(X,\iota,\rho)$
 such that $X$ is completely slope divisible. Take an arbitrary element $(X,\iota,\rho)$
 in $\mathcal{N}(\overline{\F}_p)$. 
 By \cite[Corollary 13]{MR1844205}, there exists a unique slope filtration
 $0=X_0\subset X_1\subset \cdots \subset X_s=X$.
 Put $X'=\bigoplus_{i=1}^s X_i/X_{i-1}$.
 Let us observe that $X'$ is completely slope divisible.
 Since $A_0$ and $\X$ are defined over a finite field, so is $X$.
 Therefore, by \cite[Corollary 13]{MR1844205},
 $X_i$ is also defined over a finite field. Hence \cite[Corollary 1.5]{MR1938119} tells us that
 $X'$ is completely slope divisible.
 As the category of isocrystals over $\overline{\F}_p$ is semisimple,
 the rational Dieudonn\'e modules $\mathbb{D}(X)_\Q$ and $\mathbb{D}(X')_\Q$ are canonically
 isomorphic. Thus a quasi-isogeny $f\colon X\longrightarrow X'$ is naturally induced.
 It is characterized by the property that the composite
 $X_i\hooklongrightarrow X\yrightarrow{f}X'\yrightarrow{\pr_i}X_i/X_{i-1}$ is
 the canonical surjection for each $i$.

 The $\mathcal{O}_B$-action $\iota$ on $X$ induces an $\mathcal{O}_B$-action $\iota'$ on $X'$,
 and the quasi-isogeny $f$ is $\mathcal{O}_B$-equivariant. Put $\rho'=f\circ\rho$.
 We prove that $(X',\iota',\rho')$ belongs to $\mathcal{N}(\overline{\F}_p)$.
 It suffices to see the last condition in Definition \ref{defn:N}.
 Since $(X,\iota,\rho)$ lies in $\mathcal{N}(\overline{\F}_p)$,
 we can find an element $c\in \Q_p^\times$ and an isogeny $\lambda\colon X\longrightarrow X^\vee$
 with $\height \lambda=\log_p\#(L^\vee/L)$ such that $c\lambda_0=\rho^\vee\circ\lambda\circ\rho$.
 The isogeny $\lambda\colon X\longrightarrow X^\vee$ should be compatible with the slope filtrations
 on $X$ and $X^\vee$,
 and thus it induces the isogeny $\lambda'\colon X'\longrightarrow X'^\vee$. 
 By the Dieudonn\'e theory, we can easily see that $\lambda=f^\vee\circ\lambda'\circ f$ and
 $\height\lambda'=\height\lambda=\log_p\#(L^\vee/L)$.
 Hence we have $c\lambda_0=\rho^\vee\circ\lambda\circ\rho=\rho'^\vee\circ\lambda'\circ\rho'$,
 and thus conclude that $(X',\iota',\rho')\in \mathcal{N}(\overline{\F}_p)$.

 Now, fix $(X,\iota,\rho)\in \mathcal{N}(\overline{\F}_p)$ where $X$ is completely slope divisible.
 Then, there exist an abelian variety $A_0'$ over $\overline{\F}_p$ and a $p$-quasi-isogeny
 $\phi\colon A_0\longrightarrow A_0'$ such that $\phi[p^\infty]\colon A_0[p^\infty]\longrightarrow A'_0[p^\infty]$ 
 can be identified with $\rho\colon \X\longrightarrow X$.
 By the last condition in Definition \ref{defn:N}, 
 there exists an integer $m$ such that $(\rho^\vee)^{-1}\circ p^m\lambda_0\circ \rho^{-1}$
 gives an isogeny $X\longrightarrow X^\vee$ of height $\log_p\#(L^\vee/L)$. 
 Consider the quasi-isogeny $\lambda_0'=(\phi^{\vee})^{-1}\circ p^m\lambda_0\circ \phi^{-1}\colon A_0'\longrightarrow A_0'^\vee$.
 Passing to $p$-divisible groups, we can easily observe that it is a polarization on $A_0'$.
 By construction we have $\ord_p\deg\lambda_0=\log_p\#(L^\vee/L)$.

 As $\phi$ induces $\End(A_0)\otimes_\Z\Q\yrightarrow{\cong}\End(A'_0)\otimes_\Z\Q$, $\iota_0$ induces a homomorphism
 $\iota_0'\colon \widetilde{B}\longrightarrow \End(A'_0)\otimes_\Z\Q$. Let $\widetilde{\mathcal{O}}'$ be the inverse image of $\End(A'_0)$
 under $\iota_0'$. It is an order of $\widetilde{B}$, for $\iota_0'$ is injective.
 We will show that $\widetilde{\mathcal{O}}'\otimes_{\Z}\Z_p=\mathcal{O}_B$.
 For $a\in\widetilde{B}\cap \mathcal{O}_B$, consider
 $\iota'_0(a)\in \End(A'_0)\otimes_\Z\Q$.
 Since the induced element $\iota'_0(a)[p^\infty]\in \End(A'_0[p^\infty])\otimes_{\Z_p}\Q_p$
 can be identified with $\iota(a)\in \End(X)$, it belongs to $\End(A'_0[p^\infty])$.
 Therefore we conclude that $\iota'_0(a)\in \End(A'_0)\otimes_\Z\Z_{(p)}$,
 and thus $a\in \widetilde{\mathcal{O}}'\otimes_\Z\Z_{(p)}$.
 Hence we have $(\widetilde{B}\cap \mathcal{O}_B)\otimes_{\Z_{(p)}}\Z_p\subset \widetilde{\mathcal{O}}'\otimes_\Z\Z_p$.
 On the other hand, \cite[Theorem 5.2]{MR1972204} for $R=\Z_{(p)}$ tells us that 
 $(\widetilde{B}\cap \mathcal{O}_B)\otimes_{\Z_{(p)}}\Z_p=\mathcal{O}_B$. As $\mathcal{O}_B$ is a maximal order of $B$, 
 we conclude that $\widetilde{\mathcal{O}}'\otimes_\Z\Z_p=\mathcal{O}_B$.

 Take a $\Z$-basis $e_1,\ldots,e_r$ of $\widetilde{\mathcal{O}}'$. 
 Since $\widetilde{\mathcal{O}}'\otimes_\Z\Z_p=\mathcal{O}_B$, we can find an integer $N>0$
 which is prime to $p$ such that $Ne_i^*\in \widetilde{\mathcal{O}}'$ for every $1\le i\le r$
 (recall that $\mathcal{O}_B$ is stable under $*$).
 Let $\widetilde{\mathcal{O}}$ be the $\Z$-subalgebra of $\widetilde{\mathcal{O}}'$
 generated by $Ne_i$ and $Ne_i^*$ for $1\le i\le r$.
 Then $\widetilde{\mathcal{O}}$ is an order of $\widetilde{B}$ which is contained in
 $\widetilde{\mathcal{O}}'=\iota'^{-1}_0(\End(A_0'))$, stable under $*$ and
 satisfies $\widetilde{\mathcal{O}}\otimes_\Z\Z_p=\mathcal{O}_B$.

 By construction the polarized $B$-isocrystal associated to $(A'_0,\lambda'_0,\iota'_0)$ is isomorphic to $N_b$;
 indeed the homomorphism $\mathbb{D}(A_0[p^\infty])_{\Q}\longrightarrow \mathbb{D}(A'_0[p^\infty])_{\Q}$ induced by $\phi$
 is an isomorphism of polarized $B$-isocrystals. Hence we may replace $(A_0,\lambda_0,\iota_0)$ by $(A'_0,\lambda'_0,\iota'_0)$.
\end{prf}

If $\mathcal{N}(\overline{\F}_p)=\varnothing$, then $\overline{\mathcal{N}}=\varnothing$ and Theorem \ref{thm:main} is clear. 
In the sequel, we will assume that $\mathcal{N}(\overline{\F}_p)\neq\varnothing$,
and take $(A_0,\lambda_0,\iota_0)$ and $\widetilde{\mathcal{O}}$
as in Lemma \ref{lem:replace-AV}. Put $g=\dim A_0$ and $d=\deg \lambda_0$.

\begin{defn}\label{defn:decomposition}
 For an integer $\delta$, let $\mathcal{N}^{(\delta)}$ be the open and closed formal subscheme of $\mathcal{N}$
 consisting of $(X,\iota,\rho)$ with $g^{-1}\cdot \height\rho=\delta$.
 Note that the left hand side is always an integer.
 Indeed, by the definition of $\mathcal{N}$, at least locally, 
 there exist $c\in\Q_p^\times$ and an isogeny $\lambda\colon X\longrightarrow X^\vee$ with 
 $\height \lambda=\log_p\#(L^\vee/L)$ such that $\height (c\lambda_0)=\height\rho+\height\lambda+\height\rho^\vee$.
 Since $\height\lambda_0=\log_p\#(L^\vee/L)=\height\lambda$ and $\height\rho^\vee=\height \rho$,
 we have $\height\rho=g\ord_p(c)$.
 Hence $g^{-1}\cdot \height\rho=\ord_p(c)$ is an integer.

 The formal scheme $\mathcal{N}$ is decomposed into the disjoint union $\coprod_{\delta\in\Z}\mathcal{N}^{(\delta)}$.
 Put $\overline{\mathcal{N}}{}^{(\delta)}=(\mathcal{N}^{(\delta)})^{\mathrm{red}}$.
\end{defn}

By the argument in the definition above, we can also prove the following:

\begin{lem}\label{lem:canoncal-lambda}
 Let $\delta$ be an integer. For $S\in\Nilp_{W(\overline{\F}_p)}$ and $(X,\iota,\rho)\in \mathcal{N}^{(\delta)}(S)$, 
 there exists a unique isogeny $\lambda\colon X\longrightarrow X^\vee$ with $\height\lambda=\log_p\#(L^\vee/L)$
 such that the following diagram is commutative:
 \[
 \xymatrix{%
 \X\otimes_{\overline{\F}_p}\overline{S}\ar[d]^-{p^\delta\lambda_0\otimes\id}\ar[rr]^-{\rho}&&
 X\times_S\overline{S}\ar[d]^-{\lambda\times\id}\\
 \X^\vee\otimes_{\overline{\F}_p}\overline{S}&&
 X^\vee\times_S\overline{S}\ar[ll]_-{\rho^\vee}\lefteqn{.}
 }
 \]
\end{lem}

\begin{prf}
 We have only to check that the quasi-isogeny $\lambda\colon X\longrightarrow X^\vee$ lifting
 $({\rho^\vee})^{-1}\circ (p^\delta\lambda_0\otimes\id)\circ\rho^{-1}$ is an isogeny. This is a local problem on $S$, and thus
 we may assume that there exist $c\in\Q_p^\times$ and an isogeny $\lambda'\colon X\longrightarrow X^\vee$
 with $\height\lambda'=\log_p\#(L^\vee/L)$
 such that the following diagram is commutative:
 \[
 \xymatrix{%
 \X\otimes_{\overline{\F}_p}\overline{S}\ar[d]^-{c\lambda_0\otimes\id}\ar[rr]^-{\rho}&&
 X\times_S\overline{S}\ar[d]^-{\lambda'\times\id}\\
 \X^\vee\otimes_{\overline{\F}_p}\overline{S}&&
 X^\vee\times_S\overline{S}\ar[ll]_-{\rho^\vee}\lefteqn{.}
 }
 \]
 By the argument in Definition \ref{defn:decomposition}, we have $\ord_p(c)=\delta$. 
 Therefore $c'=p^\delta c^{-1}$ lies in $\Z_p^\times$. On the other hand, we have $\lambda=c'\lambda'$, as they coincide
 over $\overline{S}$. Hence $\lambda$ is an isogeny, as desired.
\end{prf}

Next we will construct a moduli space of abelian varieties with some additional structures. 

\begin{defn}\label{defn:PEL-moduli-space}
 Fix an integer $n\ge 3$ which is prime to $p$.
 Let $\mathfrak{X}$ be the functor from the category of $\Z_p$-schemes
 to $\mathbf{Set}$
 that associates $S$ to the set of isomorphism classes of $(A,\lambda,\alpha,\iota)$ where
 \begin{itemize}
  \item $A$ is a $g$-dimensional abelian scheme over $S$,
  \item $\lambda\colon A\longrightarrow A^\vee$ is a polarization of $A$ with $\deg\lambda=d$,
  \item $\alpha\colon (\Z/n\Z)^{2g}\yrightarrow{\cong} A[n]$ is an isomorphism of group schemes
	over $S$ (an $n$-level structure on $A$),
  \item and $\iota\colon \widetilde{\mathcal{O}}\longrightarrow \End(A)$ is a homomorphism such that
	$\lambda\circ \iota(a^*)=\iota(a)^\vee\circ\lambda$ for every $a\in\widetilde{\mathcal{O}}$.
 \end{itemize}
 An isomorphism between two quadruples $(A,\lambda,\alpha,\iota)$, $(A',\lambda',\alpha',\iota')$ is an isomorphism $f\colon A\longrightarrow A'$
 such that $\lambda=f^\vee\circ\lambda'\circ f$, $\alpha'=f\circ \alpha$,
 and $\iota'(a)\circ f=f\circ\iota(a)$ for every $a\in\widetilde{\mathcal{O}}$.
\end{defn}

Note that, unlike in the definition of Shimura varieties of PEL type, we impose no compatibility condition on $\alpha$ and $\iota$.
This is to avoid an auxiliary choice of linear-algebraic data outside $p$.

\begin{prop}\label{prop:moduli-representability}
 The functor $\mathfrak{X}$ is represented by a quasi-projective scheme over $\Z_p$.
\end{prop}

\begin{prf}
 By \cite[Theorem 7.9]{MR1304906}, the moduli space of $(A,\lambda,\alpha)$ is represented by
 a quasi-projective scheme $\mathcal{A}_{g,d,n}$ over $\Z_p$.
 Furthermore, \cite[Proposition 1.3.3.7]{Kai-Wen} tells us that
 $\mathfrak{X}$ is represented by a scheme which is finite over $\mathcal{A}_{g,d,n}$.
 This concludes the proof.
\end{prf}

In the proof of Theorem \ref{thm:main}, we focus on the geometric special fiber $\overline{\mathfrak{X}}=\mathfrak{X}\otimes_{\Z_p}\overline{\F}_p$ of $\mathfrak{X}$.
For $x\in \overline{\mathfrak{X}}$ and a geometric point $\overline{x}$ lying over $x$, the corresponding quadruple 
$(A_{\overline{x}},\lambda_{\overline{x}},\alpha_{\overline{x}},\iota_{\overline{x}})$ gives
a $p$-polarized $B$-isocrystal
$\mathbb{D}(A_{\overline{x}}[p^\infty])_{\Q}$ over $\kappa(\overline{x})$.
The condition that $\mathbb{D}(A_{\overline{x}}[p^\infty])_{\Q}$ is isomorphic to $N_b\otimes_{K_0} \Frac W(\kappa(\overline{x}))$
is independent of the choice of $\overline{x}$. Indeed, the following lemma holds: 

\begin{lem}
 Let $k$ be an algebraically closed field of characteristic $p$ and $k'$ an algebraically closed extension
 field of $k$. For $p$-polarized $B$-isocrystals $N$, $N'$ over $k$, $N\cong N'$ if and only if
 $N\otimes_{\Frac W(k)}\Frac W(k')\cong N'\otimes_{\Frac W(k)}\Frac W(k')$.
\end{lem}

\begin{prf}
 Assume that there exists an isomorphism of polarized $B$-isocrystals
 \[
  f\colon N\otimes_{\Frac W(k)}\Frac W(k')\yrightarrow{\cong} N'\otimes_{\Frac W(k)}\Frac W(k').
 \] 
 By \cite[Lemma 3.9]{MR1411570}, there exists an isomorphism 
 $f'\colon N\yrightarrow{\cong}N'$ of isocrystals (without any additional structure)
 such that $f=f'\otimes \id$. Further, by \cite[Lemma 3.9]{MR1411570}
 we can observe that $f'$ is compatible with $B$-actions and polarizations
 (here we use the fact that $f$ carries the pairing on $N\otimes_{\Frac W(k)}\Frac W(k')$ to
 a $\Q_p^\times$-multiple of the pairing on $N'\otimes_{\Frac W(k)}\Frac W(k')$; see Definition \ref{defn:pol-isoc}).
 Hence $f'$ gives an isomorphism of polarized $B$-isocrystals.
\end{prf}
We write $\overline{\mathfrak{X}}{}^{(b)}$ for the set of $x\in \overline{\mathfrak{X}}$
satisfying $\mathbb{D}(A_{\overline{x}}[p^\infty])_{\Q}\cong N_b\otimes_{K_0} \Frac W(\kappa(\overline{x}))$.
The next lemma ensures that $\overline{\mathfrak{X}}{}^{(b)}$ is a locally closed
subset of $\overline{\mathfrak{X}}$.

\begin{lem}\label{lem:Newton-strata-loc-closed}
 Let $S$ be a locally noetherian scheme of characteristic $p$ and $M$ an $F$-isocrystal over $S$ 
 (\cf \cite[\S 2.1]{MR563463}, \cite[\S 3.1]{MR1411570}) endowed with a $B$-action and a non-degenerate
 alternating bilinear pairing $\langle\ ,\ \rangle\colon M\times M\longrightarrow \mathbf{1}(-1)$,
 where $\mathbf{1}(-1)$ denotes the dual of the Tate object,
 satisfying the following conditions:
 \begin{itemize}
  \item $\langle bx,y\rangle=\langle x,b^*y\rangle$ for every $b\in B$ and $x,y\in M$,
  \item and the induced map $M\otimes M\longrightarrow \mathbf{1}(-1)$ is a morphism of $F$-isocrystals.
 \end{itemize}
 Then, the set $S^{(b)}$ consisting of $s\in S$ such that
 $M_{\overline{s}}\cong N_b\otimes_{K_0}\Frac W(\kappa(\overline{s}))$
 as polarized $B$-isocrystals is locally closed in $S$.
\end{lem}

\begin{prf}
 If $\mathbf{G}$ is connected, this lemma is due to Rapoport-Richartz 
 \cite[Proposition 2.4 (iii), Theorem 3.6 (ii), Theorem 3.8]{MR1411570}. We will adapt their argument
 to our case.

 First, by Grothendieck's specialization theorem \cite[Theorem 2.3.1]{MR563463}, the set consisting of
 $s\in S$ such that $M_{\overline{s}}\cong N_b\otimes_{K_0}\Frac W(\kappa(\overline{s}))$
 as isocrystals (without any additional structure) is locally closed in $S$.
 Therefore, by replacing $S$ with this subset endowed with the induced reduced scheme structure,
 we may assume that the Newton polygons of the isocrystals $M_{\overline{s}}$ for $s\in S$ are constant.
 Further, we may also assume that $S$ is connected. In this case, we shall prove that $S^{(b)}$ is either
 $S$ or empty.
 As in the proof of \cite[Theorem 3.8]{MR1411570}, it suffices to show the following:
 \begin{quote}
  Assume that $S=\Spec k[[t]]$, where $k$ is an algebraically closed field of characteristic $p$.
  We denote by $s_1$ (resp.\ $s_0$) the generic (resp.\ closed) point of $S$.
  Then we have an isomorphism $M_{\overline{s_1}}\cong M_{s_0}\otimes_{K_0}\Frac W(\kappa(\overline{s_1}))$ of polarized $B$-isocrystals.
 \end{quote}
 Following \cite{MR1411570}, we write $R$ for the perfect closure of $k[[t]]$ and
 $a$ for the composite $\Spec R\longrightarrow \Spec k\yrightarrow{s_0}S$. 
 By \cite[Theorem 2.7.4]{MR563463} and \cite[Lemma 3.9]{MR1411570},
 there exists a unique isomorphism $f\colon M_R\yrightarrow{\cong} a^*(M)$ of $F$-isocrystals over $\Spec R$
 which induces the identity over $s_0$. We shall observe that $f$ is compatible with
 the $B$-actions and the polarizations on $M_R$ and $a^*(M)$. For $b\in B$,
 $f\circ \iota_{M_R}(b)$ and $\iota_{a^*(M)}(b)\circ f$ are elements of $\Hom(M_R,a^*(M))$, whose images
 in $\Hom(M_{s_0},M_{s_0})$ are both equal to $\iota_{M_{s_0}}(b)$. On the other hand,
 \cite[Lemma 3.9]{MR1411570} tells us that 
 the pull-back map $\Hom(M_R,a^*(M))\longrightarrow \Hom(M_{s_0},M_{s_0})$ is bijective
 (note that $M_R$ is constant, for it is isomorphic to $a^*(M)$). Thus $f$ is compatible with $B$-actions.
 The same method can be used to compare polarizations.
 Hence, by taking the fiber of $f$ at $\overline{s_1}$,
 we obtain an isomorphism $f_{\overline{s_1}}\colon M_{\overline{s_1}}\yrightarrow{\cong} M_{s_0}\otimes_{K_0}\Frac W(\kappa(\overline{s_1}))$ compatible with $B$-actions and polarizations.
 This concludes the proof.
\end{prf}
We endow $\overline{\mathfrak{X}}{}^{(b)}$ with the reduced scheme structure induced from $\overline{\mathfrak{X}}$.
By definition, the fixed $(A_0,\lambda_0,\iota_0)$ and an arbitrary $n$-level structure $\alpha_0$ on $A_0$
give an $\overline{\F}_p$-valued point of $\overline{\mathfrak{X}}{}^{(b)}$.
In the following we fix $\alpha_0$.

\subsection{Oort's leaf}\label{subsec:oort-leaf}
Here we follow the construction in \cite{MR2051612}, \cite{MR2074715} and \cite{MR2169874}.

\begin{defn}
 Let $C$ be the subset of $\overline{\mathfrak{X}}{}^{(b)}$ consisting of $x=(A_x,\lambda_x,\alpha_x,\iota_x)$ satisfying the following condition:
 \begin{itemize}
  \item for some algebraically closed field extension $k$ of $\kappa(x)$, there exists an isomorphism
	$\X\otimes_{\overline{\F}_p}k\yrightarrow{\cong}A_x[p^\infty]\otimes_{\kappa(x)}k$
	which carries $\lambda_0$ to a $\Z_p^\times$-multiple of $\lambda_x$ and $\iota_0$ to $\iota_x$.
 \end{itemize}
\end{defn}

By definition, $(A_0,\lambda_0,\alpha_0,\iota_0)$ lies in $C$. In particular $C$ is non-empty.
The goal of this subsection is to prove the following theorem, which is a generalization of \cite[Theorem 2.2]{MR2051612}.

\begin{thm}\label{thm:leaf-closed}
 The subset $C$ is closed in $\overline{\mathfrak{X}}{}^{(b)}$.
\end{thm}

The idea of our proof of Theorem \ref{thm:leaf-closed} is to use the theory of Igusa towers developed in \cite{MR2074715}.

\begin{defn}
 Let $C_{\mathrm{naive}}$ be the subset of $\overline{\mathfrak{X}}{}^{(b)}$ consisting of $x=(A_x,\lambda_x,\alpha_x,\iota_x)$
 satisfying the following condition:
 \begin{itemize}
  \item for some algebraically closed field extension $k$ of $\kappa(x)$, there exists an isomorphism
	$\X\otimes_{\overline{\F}_p}k\cong A_x[p^\infty]\otimes_{\kappa(x)}k$ as $p$-divisible groups
	(we impose no compatibility on quasi-polarizations and $\mathcal{O}_B$-actions).
 \end{itemize}
 By \cite[Theorem 2.2]{MR2051612}, $C_{\mathrm{naive}}$ is a closed subset of $\overline{\mathfrak{X}}{}^{(b)}$. We endow it with
 the induced reduced scheme structure. Let $C^\sim_{\mathrm{naive}}$ be the normalization of $C_{\mathrm{naive}}$.
\end{defn} 

Let $(\mathcal{A},\widetilde{\lambda},\widetilde{\alpha},\widetilde{\iota})$ be the universal object over $\overline{\mathfrak{X}}$
and put $\mathcal{G}=\mathcal{A}[p^\infty]$. We sometimes denote the pull-back of $\mathcal{A}$ and $\mathcal{G}$ to various schemes
by the same symbols $\mathcal{A}$ and $\mathcal{G}$.

\begin{lem}\label{lem:csd}
 The $p$-divisible group $\mathcal{G}$ over $C^\sim_{\mathrm{naive}}$ is completely slope divisible in the sense of
 \cite[Definition 1.2]{MR1938119}.
\end{lem}

\begin{prf}
 Recall that $\X$ is completely slope divisible (\cf Lemma \ref{lem:replace-AV}).
 Consider a minimal point $\eta$ of $C^\sim_{\mathrm{naive}}$.
 By the definition of $C_{\mathrm{naive}}$, there exists
 an algebraically closed extension field $k$ of $\kappa(\eta)$ such that
 $\X\otimes_{\overline{\F}_p}k\cong \mathcal{G}_\eta\otimes_{\kappa(\eta)}k$.
 Therefore \cite[Remark in p.~186]{MR1938119} tells us that $\mathcal{G}_\eta$ is completely slope divisible. Hence, by \cite[Proposition 2.3]{MR1938119}, we conclude that $\mathcal{G}$
 over $C^\sim_{\mathrm{naive}}$ is completely slope divisible.
\end{prf}

By the lemma above, we have a (unique) slope filtration $0=\mathcal{G}_0\subset \mathcal{G}_1\subset \cdots\subset \mathcal{G}_s=\mathcal{G}$.
Put $\mathcal{G}^i=\mathcal{G}_i/\mathcal{G}_{i-1}$. The action of $\mathcal{O}_B$ on $\mathcal{G}$ induces that on $\mathcal{G}^i$.
The quasi-polarization $\widetilde{\lambda}\colon \mathcal{G}\longrightarrow \mathcal{G}^\vee$ induced from the universal polarization
gives a morphism $\widetilde{\lambda}_i\colon \mathcal{G}^i\longrightarrow (\mathcal{G}^{i'})^\vee$, where $i'$ is the integer
such that $\slope \mathcal{G}^i+\slope \mathcal{G}^{i'}=1$.
Note that $\widetilde{\lambda}\colon \mathcal{G}\longrightarrow \mathcal{G}^\vee$ is an isogeny, hence so is $\widetilde{\lambda}_i$.

We can also consider $0=\X_0\subset \X_1\subset \cdots\subset \X_s=\X$, $\X^i=\X_i/\X_{i-1}$,
the $\mathcal{O}_B$-action on $\X^i$
and $\lambda_{0,i}\colon \X^i\longrightarrow (\X^{i'})^\vee$.
They appear as the fibers of the previous objects at a point of $C^\sim_{\mathrm{naive}}$ lying over 
$(A_0,\lambda_0,\alpha_0,\iota_0)\in C_{\mathrm{naive}}$.
In particular, $\lambda_{0,i}$ is an isogeny.

\begin{defn}\label{defn:naive-Igusa}
 For an integer $m\ge 0$, let $C^\sim_{\mathrm{naive},m}$ be the functor from the category of
 $C^\sim_{\mathrm{naive}}$-schemes to $\mathbf{Set}$
 that associates $S$ to the set of $\{j_{m,i}\}_{1\le i\le s}$
 where $j_{m,i}\colon \X^i[p^m]\otimes_{\overline{\F}_p}S\yrightarrow{\cong} \mathcal{G}^i[p^m]\times_{C^\sim_{\mathrm{naive}}}S$
 is an isomorphism of group schemes over $S$.
 This functor is represented by a scheme of finite type over $C^\sim_{\mathrm{naive}}$. 
 We denote the universal isomorphisms over $C^\sim_{\mathrm{naive},m}$ by $j^{\mathrm{univ}}_{m,i}$.

 Let $C^{\sim\prime}_{\mathrm{naive},m}$ be the intersection of the scheme-theoretic images of
 $C^\sim_{\mathrm{naive},m'}\longrightarrow C^\sim_{\mathrm{naive},m}$ for $m'\ge m$ and put
 $\Ig^\sim_{\mathrm{naive},m}=(C^{\sim\prime}_{\mathrm{naive},m})^{\mathrm{red}}$.
 We also denote the restriction of $j^{\mathrm{univ}}_{m,i}$ to $\Ig^\sim_{\mathrm{naive},m}$ by the same symbol $j^{\mathrm{univ}}_{m,i}$.
\end{defn}

The following proposition is due to Harris-Taylor and Mantovan:

\begin{prop}\label{prop:Igusa-fin-etale}
 For each $m\ge 0$, the morphism $\Ig^\sim_{\mathrm{naive},m+1}\longrightarrow \Ig^\sim_{\mathrm{naive},m}$ is finite \'etale and surjective.
\end{prop}

\begin{prf}
 This is essentially proved in \cite[Proposition 3.3]{MR2074715}, but we should be careful since
 \cite[Proposition 3.3]{MR2074715} is stated in the case of unitary Shimura varieties.
 Recall that its proof is a combination of \cite[Proposition II.1.7]{MR1876802}
 and \cite[Lemma 3.4]{MR2074715}. The former is valid for any reduced excellent scheme over $\overline{\F}_p$.
 The latter requires the fact that the completed local ring
 $\mathcal{O}^{\wedge}_{C^\sim_{\mathrm{naive}},x}$ is a normal integral domain
 for each closed point $x\in C^\sim_{\mathrm{naive}}$.
 This is true because $C^\sim_{\mathrm{naive}}$ is normal and excellent
 (\cf \cite[IV, 7.8.3 (v)]{EGA}).
\end{prf}

\begin{defn}\label{defn:refined-Igusa}
 For $m\ge 0$ and a $C^\sim_{\mathrm{naive}}$-scheme $S$, let $C^\sim_m(S)\subset \Ig^\sim_{\mathrm{naive},m}(S)$
 be the subset consisting of $\{j_{m,i}\}\in \Ig^\sim_{\mathrm{naive},m}(S)$ satisfying
 the following two conditions:
 \begin{itemize}
  \item The isomorphisms $\{j_{m,i}\}$ preserve the $(\Z/p^m\Z)^\times$-homothety
	classes of the quasi-polarizations. More precisely, for each $i$, the following diagram commutes
	up to multiplication by $(\Z/p^m\Z)^\times$:
	\[
	\xymatrix{%
	\X^i[p^m]\otimes_{\overline{\F}_p}S\ar[rr]_-{\cong}^-{j_{m,i}}\ar[d]^-{\lambda_{0,i}}&& \mathcal{G}^i[p^m]\times_{C^\sim_{\mathrm{naive}}}S\ar[d]^-{\widetilde{\lambda}_i}\\
	(\X^{i'})^\vee[p^m]\otimes_{\overline{\F}_p}S&& (\mathcal{G}^{i'})^{\vee}[p^m]\times_{C^\sim_{\mathrm{naive}}}S\ar[ll]^-{\cong}_-{(j_{m,i'})^\vee}\lefteqn{.}
	}
	\]
	Here $i'$ is the integer such that 
	$\slope \mathcal{G}^i+\slope \mathcal{G}^{i'}=\slope \X^i+\slope \X^{i'}=1$.
	Recall that $\widetilde{\lambda}_i\colon \mathcal{G}^i\longrightarrow (\mathcal{G}^{i'})^\vee$ and 
	$\lambda_{0,i}\colon \X^i\longrightarrow (\X^{i'})^\vee$ are isogenies
	(see the explanation before Definition \ref{defn:naive-Igusa}), so that
	the vertical arrows are induced.
  \item For each $i$, the isomorphism $j_{m,i}\colon \X^i[p^m]\otimes_{\overline{\F}_p}S\yrightarrow{\cong} \mathcal{G}^i[p^m]\times_{C^\sim_{\mathrm{naive}}}S$ preserves the $\mathcal{O}_B$-actions
	explained before Definition \ref{defn:naive-Igusa}.
 \end{itemize}
 Clearly $S\longmapsto C^\sim_m(S)$ is a functor represented by a closed subscheme of
 $\Ig^\sim_{\mathrm{naive},m}$.
 Moreover, let $C^{\sim\prime}_m$ be the intersection of the scheme-theoretic images of $C^\sim_{m'}\longrightarrow C^\sim_m$ for $m'\ge m$.
 Put $\Ig^\sim_m=(C^{\sim\prime}_m)^{\mathrm{red}}$ and $C^\sim=\Ig^\sim_0$.
\end{defn}

\begin{cor}\label{cor:refined-Igusa-fin-surj}
 For each $m\ge 0$, the morphism $\Ig^\sim_{m+1}\longrightarrow \Ig^\sim_m$ is finite and surjective.
\end{cor}

\begin{prf}
 Since $\Ig^\sim_m$ (resp.\ $\Ig^\sim_{m+1}$) is a closed subscheme of $\Ig^\sim_{\mathrm{naive},m}$ (resp.\ $\Ig^\sim_{\mathrm{naive},m+1}$), 
 the finiteness follows from Proposition \ref{prop:Igusa-fin-etale}. To show the surjectivity, note that
 there exists an integer $m'\ge m+1$ such that $C^{\sim\prime}_{m+1}$ coincides with
 the scheme-theoretic image of
 $C^\sim_{m'}\longrightarrow C^\sim_{m+1}$, for $C^\sim_{m+1}$ is a noetherian scheme.
 Take $x\in C^{\sim\prime}_m$. As $C^\sim_{m'}\longrightarrow C^\sim_m$ is finite by Proposition \ref{prop:Igusa-fin-etale},
 $x$ lies in the set-theoretic image of this map. Hence we can find $x'$ in $C^\sim_{m'}$ which is mapped to $x$.
 The image of $x'$ in $C^\sim_{m+1}$ lies in $C^{\sim\prime}_{m+1}$, and is mapped to $x$.
 Therefore $\Ig^\sim_{m+1}=C^{\sim\prime}_{m+1}\longrightarrow C^{\sim\prime}_m=\Ig^\sim_m$ is surjective, as desired.
\end{prf}

The scheme $C^\sim$ is a closed subscheme of $C^\sim_{\mathrm{naive}}$. Therefore,
to prove Theorem \ref{thm:leaf-closed}, it suffices to show the following proposition.

\begin{prop}\label{prop:leaf-image}
 The image of $C^\sim$ under the finite morphism $C^\sim_{\mathrm{naive}}\longrightarrow C_{\mathrm{naive}}$ coincides with $C$.
\end{prop}

\begin{prf}
 Take $x\in C^\sim_{\mathrm{naive}}$. First assume that the image of $x$ in $C_{\mathrm{naive}}$ lies in $C$.
 Then there exist an algebraically closed extension field $k$ of $\kappa(x)$ and an isomorphism
 $j\colon \X\otimes_{\overline{\F}_p}k\yrightarrow{\cong} \mathcal{G}_x\otimes_{\kappa(x)}k$ preserving
 the $\Z_p^\times$-homothety classes of the quasi-polarizations
 and the $\mathcal{O}_B$-actions. It induces an isomorphism
 $j_{m,i}\colon \X^i[p^m]\otimes_{\overline{\F}_p}k\yrightarrow{\cong} \mathcal{G}^i_x[p^m]\otimes_{\kappa(x)}k$
 for each $m\ge 0$ and $1\le i\le s$.
 Therefore we obtain a system of maps $\{\Spec k\longrightarrow C^\sim_{\mathrm{naive},m}\}_{m\ge 0}$ compatible with projections.
 It is easy to check that the morphism $\Spec k\longrightarrow C^\sim_{\mathrm{naive},m}$ factors through $\Ig^\sim_m$.
 In particular, the image of $\Spec k\longrightarrow C^\sim_{\mathrm{naive}}$, which is nothing but $x$, lies in $C^\sim$.

 Conversely assume that $x$ lies in $C^\sim$. By Corollary \ref{cor:refined-Igusa-fin-surj}, we can take a system of points
 $\{x_m\in \Ig^\sim_m\}_{m\ge 0}$ compatible with projections such that $x_0=x$.
 Let $k$ be an algebraic closure of $\varinjlim_m\kappa(x_m)$. Then, for each $m$ we have a collection of isomorphisms
 $\{j_{m,i}\colon \X^i[p^m]\otimes_{\overline{\F}_p}k\yrightarrow{\cong} \mathcal{G}^i_x[p^m]\otimes_{\kappa(x)}k\}_{1\le i\le s}$
 compatible with the change of $m$.
 By \cite[Corollary 11]{MR1844205}, the slope filtrations on $\X$ and $\mathcal{G}_x\otimes_{\kappa(x)}k$ split canonically.
 Namely, we have isomorphisms $\X\cong \bigoplus_{i=1}^s\X^i$ and $\mathcal{G}_x\otimes_{\kappa(x)}k\cong \bigoplus_{i=1}^s\mathcal{G}^i_x\otimes_{\kappa(x)}k$ which are compatible with the quasi-polarizations and the $\mathcal{O}_B$-actions.
 Hence, $\{j_{m,i}\}_{1\le i\le s}$ induces an isomorphism 
 $j_m\colon \X[p^m]\otimes_{\overline{\F}_p}k\yrightarrow{\cong} \mathcal{G}_x[p^m]\otimes_{\kappa(x)}k$
 compatible with the $(\Z/p^m\Z)^\times$-homothety classes of the quasi-polarizations,
 the $\mathcal{O}_B$-actions and the change of $m$.
 By taking inductive limit with respect to $m$, we obtain an isomorphism $j\colon \X\otimes_{\overline{\F}_p}k\yrightarrow{\cong}\mathcal{G}_x\otimes_{\kappa(x)}k$ compatible with the $\mathcal{O}_B$-actions.
 Since the projective limit of a projective system consisting of finite sets is non-empty,
 $j$ preserves the $\Z_p^\times$-homothety classes of the quasi-polarizations.
 This means that the image of $x$ in $C_{\mathrm{naive}}$ belongs to $C$.

 Noting that $C^\sim_{\mathrm{naive}}\longrightarrow C_{\mathrm{naive}}$ is surjective, 
 we conclude the proof.
\end{prf}

Now the proof of Theorem \ref{thm:leaf-closed} is complete.
We endow the closed subset $C\subset \overline{\mathfrak{X}}{}^{(b)}$ with the induced reduced scheme structure.

\begin{prop}
 The scheme $C$ is smooth over $\overline{\F}_p$.
\end{prop}

\begin{prf}
 It can be proved in the same way as \cite[Theorem 3.13 (i)]{MR2051612}.
\end{prf}

By this proposition, we can apply the same construction as in Definition \ref{defn:naive-Igusa}
and Definition \ref{defn:refined-Igusa} to $\mathcal{G}/C$.
Hence we first obtain the naive Igusa tower $\{\Ig_{\mathrm{naive},m}\}$, and after modifying it, the Igusa tower $\{\Ig_m\}$. The following proposition can be proved in the same way as Corollary
\ref{cor:refined-Igusa-fin-surj} and Proposition \ref{prop:leaf-image}.

\begin{prop}\label{prop:genuine-Igusa-tower}
 \begin{enumerate}
  \item The transition maps of $\{\Ig_m\}$ are finite and surjective.
  \item We have $\Ig_0=C$.
 \end{enumerate}
\end{prop}

\subsection{Almost product structure}
In this subsection, an element $(X,\iota,\rho)$ of $\mathcal{N}(S)$ for $S\in\Nilp_{W(\overline{\F}_p)}$
will be denoted by $(X,\rho)$.
For integers $r$, $m$ with $m\ge 0$, let $\mathcal{N}^{r,m}$ be the closed formal subscheme of $\mathcal{N}$ consisting of
$(X,\rho)$ such that $p^r\rho$ is an isogeny and $\Ker (p^r\rho)$ is killed by $p^m$.
We put $\overline{\mathcal{N}}{}^{r,m}=(\mathcal{N}^{r,m})^{\mathrm{red}}$.
It is a scheme of finite type over $\overline{\F}_p$ (\cf \cite[Corollary 2.31]{MR1393439}).

As in \cite{MR2051612}, \cite{MR2074715} and \cite{MR2169874}, we will construct a morphism 
\[
 \Ig_m\times_{\overline{\F}_p}\overline{\mathcal{N}}{}^{r,m}\longrightarrow \overline{\mathfrak{X}}{}^{(b)}.
\]
Recall that we denote by $(\mathcal{A},\widetilde{\lambda},\widetilde{\alpha},\widetilde{\iota})$ the universal object over $\Ig_m$.
For an integer $N\ge 0$, let $\mathcal{A}^{(p^N)}$ be the pull-back of $\mathcal{A}$ by the $N$th power of the absolute Frobenius morphism
$\Fr\colon \Ig_m\longrightarrow \Ig_m$.
By \cite[Lemma 4.1]{MR2074715}, there exists an integer $\delta_m\ge 0$
depending on $m$ such that for $N\ge \delta_m$ we have a canonical isomorphism
\[
 \mathcal{A}^{(p^N)}[p^m]=\mathcal{G}^{(p^N)}[p^m]\cong \bigoplus_{i=1}^s(\mathcal{G}^i)^{(p^N)}[p^m].
\]
Then, the universal Igusa structure $\{j_{m,i}^{\mathrm{univ}}\}$ over $\Ig_m$ gives an isomorphism
\[
 \mathcal{A}^{(p^N)}[p^m]\cong \bigoplus_{i=1}^s(\mathcal{G}^i)^{(p^N)}[p^m]\yleftarrow[\cong]{\bigoplus_i j_{m,i}^{\mathrm{univ}}}
 \bigoplus_{i=1}^s(\X^i)^{(p^N)}[p^m]\otimes_{\overline{\F}_p}\Ig_m\cong \X^{(p^N)}[p^m]\otimes_{\overline{\F}_p}\Ig_m,
\]
where $(\X^i)^{(p^N)}$ and $\X^{(p^N)}$ are the pull-back of $\X^i$ and $\X$ by the $N$th power of
the absolute Frobenius morphism on $\Spec \overline{\F}_p$.
It preserves the $(\Z/p^m\Z)^\times$-homothety classes of the quasi-polarizations
and the $\mathcal{O}_B$-actions.

On the other hand, over $\overline{\mathcal{N}}{}^{r,m}$ we have the universal $p$-divisible group $\widetilde{X}$ with an $\mathcal{O}_B$-action
and the universal $\mathcal{O}_B$-quasi-isogeny $\widetilde{\rho}\colon \X\otimes_{\overline{\F}_p}\overline{\mathcal{N}}{}^{r,m}\longrightarrow \widetilde{X}$. By the definition of $\overline{\mathcal{N}}{}^{r,m}$, $p^r\widetilde{\rho}$ is an isogeny and its kernel
$\Ker (p^r\widetilde{\rho})$ is contained in $\X[p^m]\otimes_{\overline{\F}_p}\overline{\mathcal{N}}{}^{r,m}$.
Hence $\Ker (p^r\widetilde{\rho})^{(p^N)}$ is a finite flat subgroup scheme of
$\X^{(p^N)}[p^m]\otimes_{\overline{\F}_p}\overline{\mathcal{N}}{}^{r,m}$.

Now consider the abelian scheme $\pr_1^*\mathcal{A}^{(p^N)}$ on $\Ig_m\times_{\overline{\F}_p}\overline{\mathcal{N}}{}^{r,m}$.
Under the isomorphism $\pr_1^*\mathcal{A}^{(p^N)}[p^m]\cong \X^{(p^N)}[p^m]\otimes_{\overline{\F}_p}(\Ig_m\times_{\overline{\F}_p}\overline{\mathcal{N}}{}^{r,m})$, $\pr_2^*\Ker (p^r\widetilde{\rho})^{(p^N)}$ corresponds to a finite flat group scheme $\mathcal{H}$ of
$\pr_1^*\mathcal{A}^{(p^N)}[p^m]$. We put $\mathcal{A}'=(\pr_1^*\mathcal{A}^{(p^N)})/\mathcal{H}$.
It is an abelian scheme over $\Ig_m\times_{\overline{\F}_p}\overline{\mathcal{N}}{}^{r,m}$ endowed with a $p$-isogeny
$\phi\colon \pr_1^*\mathcal{A}^{(p^N)}\longrightarrow \mathcal{A}'$. 
We will find additional structures on $\mathcal{A}'$ so that they give an element of $\overline{\mathfrak{X}}(\Ig_m\times_{\overline{\F}_p}\overline{\mathcal{N}}{}^{r,m})$.
The $\widetilde{\mathcal{O}}$-action on $\pr_1^*\mathcal{A}^{(p^N)}$ induced from $\widetilde{\iota}$ gives an $\widetilde{\mathcal{O}}$-action
$\iota'$ on $\mathcal{A}'$ by the isogeny $\phi$.
Indeed, it suffices to observe that $\mathcal{H}$ is stable under the $\widetilde{\mathcal{O}}$-action,
which follows from the fact that $\Ker (p^r\widetilde{\rho})$ is stable under the $\mathcal{O}_B$-action.
As $n$ is assumed to be prime to $p$, the quasi-isogeny $p^{-r}\phi$ induces an isomorphism 
$\pr_1^*\mathcal{A}^{(p^N)}[n]\longrightarrow \mathcal{A}'[n]$ on $n$-torsion points.
Let $\alpha'$ be the level structure on $\mathcal{A}'$ induced
from $\widetilde{\alpha}$ by this isomorphism. 
We shall construct a polarization $\lambda'$ on $\mathcal{A}'$.
Recall that we have a decomposition $\overline{\mathcal{N}}=\coprod_{\delta\in\Z}\overline{\mathcal{N}}{}^{(\delta)}$
into open and closed subschemes.
Put $\overline{\mathcal{N}}{}^{(\delta),r,m}=\overline{\mathcal{N}}{}^{(\delta)}\cap \overline{\mathcal{N}}{}^{r,m}$.
It suffices to construct $\lambda'$ over $\Ig_m\times_{\overline{\F}_p}\overline{\mathcal{N}}{}^{(\delta),r,m}$ for each $\delta\in\Z$.
By Lemma \ref{lem:canoncal-lambda}, on $\overline{\mathcal{N}}{}^{(\delta),r,m}$ there exists an isogeny
$\lambda_{\widetilde{X}}\colon \widetilde{X}\longrightarrow \widetilde{X}^\vee$ with height $\log_p\#(L^\vee/L)$
such that the following diagram is commutative:
\[
 \xymatrix{%
 \X\otimes_{\overline{\F}_p}\overline{\mathcal{N}}{}^{r,m}\ar[d]^-{p^\delta\lambda_0\otimes\id}\ar[rr]^-{\widetilde{\rho}}&&
 \widetilde{X}\ar[d]^-{\lambda_{\widetilde{X}}}\\
 \X^\vee\otimes_{\overline{\F}_p}\overline{\mathcal{N}}{}^{r,m}&&
 \widetilde{X}^\vee\ar[ll]_-{\widetilde{\rho}^\vee}\lefteqn{.}
 }
\]
Put $\lambda'=(\phi^\vee)^{-1}\circ p^{\delta+2r}\pr_1^*\widetilde{\lambda}^{(p^N)}\circ\phi^{-1}$.
Let us observe that the quasi-isogeny $\lambda'[p^\infty]\colon \mathcal{A}'[p^\infty]\longrightarrow \mathcal{A}'^\vee[p^\infty]$ induced by $\lambda'$ is
an isogeny with height $\log_p\#(L^\vee/L)$. 
As $\Ig_m\times_{\overline{\F}_p}\overline{\mathcal{N}}{}^{r,m}$ is reduced, 
by \cite[Proposition 2.9]{MR1393439}, it suffices to show that, for each point $z=(x,y)$
in $(\Ig_m\times_{\overline{\F}_p}\overline{\mathcal{N}}{}^{r,m})(\overline{\F}_p)$, the quasi-isogeny
$\lambda'_z[p^\infty]$ is an isogeny with height $\log_p\#(L^\vee/L)$.
By \cite[Corollary 11]{MR1844205}, we have an isomorphism
$\mathcal{A}_x^{(p^N)}[p^\infty]=\mathcal{G}_x^{(p^N)}\cong \bigoplus_{i=1}^s(\mathcal{G}_x^i)^{(p^N)}$.
Hence, Proposition \ref{prop:genuine-Igusa-tower} i) ensures that the isomorphism
$\mathcal{A}^{(p^N)}_x[p^m]\cong \X^{(p^N)}[p^m]$, the specialization at $x$
of the isomorphism used above, can be extended to an isomorphism
$\mathcal{A}^{(p^N)}_x[p^\infty]\cong \X^{(p^N)}$ between $p$-divisible groups
which preserves the $\Z_p^\times$-homothety classes of the quasi-polarizations
 and the $\mathcal{O}_B$-actions
 (\cf the proof of Proposition \ref{prop:leaf-image}).
Under this isomorphism, $\mathcal{A}'_z[p^\infty]$
is identified with $\widetilde{X}^{(p^N)}_y$, and $\lambda'_z[p^\infty]$ fits into the following
diagram, which is commutative up to $\Z_p^\times$-multiplication:
\[
 \xymatrix{%
 \X^{(p^N)}\ar[rr]^-{\phi_z[p^\infty]=p^r\widetilde{\rho}_y^{(N)}}\ar[d]^-{p^{\delta+2r}\lambda_0^{(p^N)}}&\hspace{50pt}& \widetilde{X}^{(p^N)}_y\ar[d]^-{\lambda'_z[p^\infty]}\\
 (\X^\vee)^{(p^N)}&\hspace{50pt}& (\widetilde{X}_y^\vee)^{(p^N)}\ar[ll]_-{\phi_z^\vee[p^\infty]=p^r(\widetilde{\rho}_y^\vee)^{(N)}}\lefteqn{.}
 }
\]
Hence $\lambda'_z[p^\infty]$ is identified with a $\Z_p^\times$-multiple of
$\lambda^{(p^N)}_{\widetilde{X},y}$, which is an isogeny
with height $\log_p\#(L^\vee/L)$.
As $\deg\phi$ is a power of $p$, we conclude that $\lambda'$ is an isogeny with degree $d$.
Thus it gives a desired polarization on $\mathcal{A}'$.
Obviously the quadruple $(\mathcal{A}',\lambda',\alpha',\iota')$ belongs to 
$\overline{\mathfrak{X}}(\Ig_m\times_{\overline{\F}_p}\overline{\mathcal{N}}{}^{r,m})$.

\begin{defn}
 Let $\pi_N\colon \Ig_m\times_{\overline{\F}_p}\overline{\mathcal{N}}{}^{r,m}\longrightarrow \overline{\mathfrak{X}}$ be the morphism 
 determined by the quadruple $(\mathcal{A}',\lambda',\alpha',\iota')$. 
 This morphism factors through $\overline{\mathfrak{X}}{}^{(b)}$
 (note that the $N$th power of the relative Frobenius morphism gives an isogeny between
 $\X$ and $\X^{(p^N)}$).
\end{defn}

The following lemma, which is an analogue of \cite[Proposition 4.3]{MR2074715}, is easily verified.

\begin{lem}\label{lem:pi_N-property}
 \begin{enumerate}
  \item For an integer $N\ge \delta_m$, we have $\pi_{N+1}=\Frob_p\circ \pi_N$,
	where $\Frob_p\colon \overline{\mathfrak{X}}\longrightarrow \overline{\mathfrak{X}}$ is
	the $p$th power Frobenius morphism over $\overline{\F}_p$ (namely, the base change to $\overline{\F}_p$ of the absolute Frobenius
	morphism on $\mathfrak{X}\otimes_{\Z_p}\F_p$).
  \item For an integer $N\ge \max\{\delta_m,\delta_{m+1}\}$, the following diagram is commutative:
	\[
	 \xymatrix{%
	\Ig_{m+1}\times_{\overline{\F}_p}\overline{\mathcal{N}}{}^{r,m}\ar[r]^-{\subset}\ar[d]&
	\Ig_{m+1}\times_{\overline{\F}_p}\overline{\mathcal{N}}{}^{r,m+1}\ar[r]^-{\pi_N}&\overline{\mathfrak{X}}{}^{(b)}\ar@{=}[d]\\
	\Ig_m\times_{\overline{\F}_p}\overline{\mathcal{N}}{}^{r,m}\ar[rr]^-{\pi_N}&& \overline{\mathfrak{X}}{}^{(b)}\lefteqn{.}
	}
	\]
 \end{enumerate}
\end{lem}

Let $k$ be an algebraically closed field containing $\overline{\F}_p$. At the level of $k$-valued points, we can define
a variant $\Pi\colon \Ig_m(k)\times\overline{\mathcal{N}}{}^{r,m}(k)\longrightarrow \overline{\mathfrak{X}}{}^{(b)}(k)$ of $\pi_N$.
Let $x=(A,\lambda,\alpha,\iota,\{j_{m,i}\}_{1\le i\le s})$ be an element of $\Ig_m(k)$ and $y=(X,\rho)$ be an element of
$\overline{\mathcal{N}}{}^{r,m}(k)$. In this case, by \cite[Corollary 11]{MR1844205}, we have a canonical isomorphism
$A[p^\infty]=\mathcal{G}_x\cong \bigoplus_{i=1}^s\mathcal{G}^i_x$. Hence $\{j_{m,i}\}_{1\le i\le s}$ induces an isomorphism
$A[p^m]\cong \X[p^m]\otimes_{\overline{\F}_p}k$ preserving the $(\Z/p^m\Z)^\times$-homothety classes of the quasi-polarizations and the $\mathcal{O}_B$-actions
(in this argument, we need no restriction on $m\ge 0$). By this isomorphism, the finite subgroup 
scheme $\Ker (p^r\rho)$ of
$\X[p^m]\otimes_{\overline{\F}_p}k$ corresponds to a subgroup scheme $H$ of $A[p^m]$. Put $A'=A/H$.
In the same way as in the definition of $\pi_N$, we can find a polarization $\lambda'$, a level structure $\alpha'$ and
an $\widetilde{\mathcal{O}}$-action $\iota'$ on $A'$ induced from $\lambda$, $\alpha$ and $\iota$ respectively,
so that $(A',\lambda',\alpha',\iota')$ belongs to the set $\overline{\mathfrak{X}}{}^{(b)}(k)$.
We define the map $\Pi\colon \Ig_m(k)\times\overline{\mathcal{N}}{}^{r,m}(k)\longrightarrow \overline{\mathfrak{X}}{}^{(b)}(k)$
by $\Pi(x,y)=(A',\lambda',\alpha',\iota')$.
The following lemma is clear.

\begin{lem}\label{lem:pi-Pi}
 \begin{enumerate}
  \item $\Pi\colon \Ig_m(k)\times\overline{\mathcal{N}}{}^{r,m}(k)\longrightarrow \overline{\mathfrak{X}}{}^{(b)}(k)$ is compatible with
	the change of $m$.
  \item $\pi_N=\Frob_p^N\circ \Pi$.
 \end{enumerate}
\end{lem}

\begin{lem}\label{lem:pointwise-surjection}
 Let $k$ be an algebraically closed field containing $\overline{\F}_p$.
 For every point $x=(A,\lambda,\alpha,\iota)$ of $\overline{\mathfrak{X}}{}^{(b)}(k)$, there exists an integer $m\ge 0$ such that
 $x$ is contained in the image of $\Pi\colon \Ig_m(k)\times\overline{\mathcal{N}}{}^{0,m}(k)\longrightarrow \overline{\mathfrak{X}}{}^{(b)}(k)$.
\end{lem}

\begin{prf}
 As $x$ lies in $\overline{\mathfrak{X}}{}^{(b)}$, there is a quasi-isogeny 
 $\rho\colon \X\otimes_{\overline{\F}_p}k\longrightarrow A[p^\infty]$
 which is compatible with the $\mathcal{O}_B$-actions
 and preserves the quasi-polarizations up to $\Q_p^\times$-multiplication.
 Replacing $\rho$ by $p^\nu\rho$ if necessary, we may assume that $\rho$ is an isogeny.
 Take an integer $m\ge 0$ such that $\Ker\rho$ is killed by $p^m$.
 Then, there exists an isogeny $\xi\colon A[p^\infty]\longrightarrow \X\otimes_{\overline{\F}_p}k$ with
 $\xi\circ\rho=p^m$.
 Put $A'=A/\Ker\xi$.
 Then $A'[p^\infty]$ can be identified with $\X\otimes_{\overline{\F}_p}k$, and
 there exists a $p$-isogeny $\phi\colon A'\longrightarrow A$ corresponding to $\rho$.
 As in the construction of $\pi_N$, we can observe that the additional structures $\lambda$, $\alpha$, $\iota$ on $A$
 induce additional structures $\lambda'$, $\alpha'$, $\iota'$ on $A'$ by the isogeny $\phi$ so that
 $(A',\lambda',\alpha',\iota')$ gives an element of $C(k)$.
 Let $\{j_{m,i}\}_{1\le i\le m}$ be the Igusa structure on $A'[p^m]$ that comes from the identification 
 $A'[p^\infty]=\X\otimes_{\overline{\F}_p}k$, and $y$ be the point $(A',\lambda',\alpha',\iota',\{j_{m,i}\})$ in $\Ig_m(k)$.
 Now it is easy to check that $z=(A[p^\infty],\rho)$ lies in $\overline{\mathcal{N}}{}^{0,m}(k)$ and $\Pi(y,z)=x$. 
\end{prf}

\begin{prop}\label{prop:pi-surjection}
 \begin{enumerate}
  \item There exist integers $m\ge 0$ and $N\ge \delta_m$ such that the morphism 
	$\pi_N\colon \Ig_m\times_{\overline{\F}_p}\overline{\mathcal{N}}{}^{0,m}\longrightarrow \overline{\mathfrak{X}}{}^{(b)}$ is surjective.
  \item For $m$ as in i) and an algebraically closed field $k$ containing $\overline{\F}_p$,
	the map $\Pi\colon \Ig_m(k)\times\overline{\mathcal{N}}{}^{0,m}(k)\longrightarrow \overline{\mathfrak{X}}{}^{(b)}(k)$ is surjective.
 \end{enumerate}
\end{prop}

\begin{prf}
 For each $m\ge 0$, take an integer $N_m\ge \delta_m$ so that $N_0<N_1<N_2<\cdots$, and put
 $T_m=(\Frob_p^{N_m})^{-1}\pi_{N_m}(\Ig_m\times_{\overline{\F}_p}\overline{\mathcal{N}}{}^{0,m})$.
 By Lemma \ref{lem:pi_N-property}, we have
 \begin{align*}
  T_m&\subset (\Frob_p^{N_{m+1}})^{-1}\pi_{N_{m+1}}(\Ig_m\times_{\overline{\F}_p}\overline{\mathcal{N}}{}^{0,m})
 =(\Frob_p^{N_{m+1}})^{-1}\pi_{N_{m+1}}(\Ig_{m+1}\times_{\overline{\F}_p}\overline{\mathcal{N}}{}^{0,m})\\
  &\subset (\Frob_p^{N_{m+1}})^{-1}\pi_{N_{m+1}}(\Ig_{m+1}\times_{\overline{\F}_p}\overline{\mathcal{N}}{}^{0,m+1})=T_{m+1}.
 \end{align*}
 Since $\overline{\mathfrak{X}}{}^{(b)}$ is a scheme of finite type over $\overline{\F}_p$, its underlying topological space is a spectral space
 in the sense of \cite[\S 0]{MR0251026}. As in \cite[\S 2]{MR0251026}, we consider the patch topology on $\overline{\mathfrak{X}}{}^{(b)}$,
 under which $\overline{\mathfrak{X}}{}^{(b)}$ becomes a compact space (\cite[Theorem 1]{MR0251026}).
 As $\overline{\mathcal{N}}{}^{0,m}$ is a scheme of finite type over $\overline{\F}_p$, $T_m$ is a constructible subset of $\overline{\mathfrak{X}}{}^{(b)}$. In particular, it is an open set of $\overline{\mathfrak{X}}{}^{(b)}$ with respect to the patch topology.
 On the other hand, Lemma \ref{lem:pi-Pi} ii) and Lemma \ref{lem:pointwise-surjection} tell us that
 $\bigcup_{m=0}^\infty T_m=\overline{\mathfrak{X}}{}^{(b)}$. Hence there exists an integer $m\ge 0$ such that
 $T_m=\overline{\mathfrak{X}}{}^{(b)}$. As $\Frob_p$ is surjective, 
 we have $\pi_{N_m}(\Ig_m\times_{\overline{\F}_p}\overline{\mathcal{N}}{}^{0,m})=\overline{\mathfrak{X}}{}^{(b)}$.
 This concludes the proof of i).

 ii) follows from i), Lemma \ref{lem:pi-Pi} ii) and the injectivity of $\Frob_p$ at the level of $k$-valued points.
\end{prf}

\begin{prop}\label{prop:fiber-orbit}
 Let $m,m'\ge 0$ be integers, and consider elements 
 $(x,y)\in \Ig_m(\overline{\F}_p)\times\overline{\mathcal{N}}{}^{0,m}(\overline{\F}_p)$ and 
 $(x',y')\in \Ig_{m'}(\overline{\F}_p)\times\overline{\mathcal{N}}{}^{0,m'}(\overline{\F}_p)$.
 If $\Pi(x,y)=\Pi(x',y')$, then there exists $h\in J$ such that $y'=hy$.
\end{prop}

\begin{prf}
 By Proposition \ref{prop:genuine-Igusa-tower} i), we may assume that $m=m'$. 
 Let $(A,\lambda,\alpha,\iota,\{j_{m,i}\})$, $(A',\lambda',\alpha',\iota',\{j'_{m,i}\})$,
 $(X,\rho)$, $(X',\rho')$ be the objects corresponding to $x$, $x'$, $y$, $y'$, respectively.
 Recall that $\{j_{m,i}\}$ induces an isomorphism $j_m\colon \X[p^m]\yrightarrow{\cong}A[p^m]$.
 Again by Proposition \ref{prop:genuine-Igusa-tower} i),
 we can extend $j_m$ to an isomorphism $j\colon \X\yrightarrow{\cong}A[p^\infty]$
 which preserves the $\Z_p^\times$-homothety classes of the quasi-polarizations
 and the $\mathcal{O}_B$-actions
 (\cf the proof of Proposition \ref{prop:leaf-image}).
 Similarly we have isomorphisms $j'_m\colon \X[p^m]\yrightarrow{\cong}A'[p^m]$ and $j'\colon \X\yrightarrow{\cong}A'[p^\infty]$.
 As $\Pi(x,y)=\Pi(x',y')$, there exists an isomorphism $A/j(\Ker\rho)\cong A'/j'(\Ker\rho')$ compatible with
 the additional structures. On the other hand, $j$ (resp.\ $j'$) induces an isomorphism $X\cong \X/\Ker\rho\yrightarrow{\cong} (A/j(\Ker\rho))[p^\infty]$
 (resp.\ $X'\cong\X/\Ker\rho'\yrightarrow{\cong} (A'/j'(\Ker\rho'))[p^\infty]$).
 Hence we obtain an isomorphism $f\colon X\yrightarrow{\cong}X'$, which is compatible with the $\mathcal{O}_B$-actions.
 Consider the quasi-isogeny $h=\rho'^{-1}\circ f\circ\rho\colon \X\longrightarrow \X$. It is straightforward to check that
 $h$ in fact gives an element of $J$ (use the fact that the isogenies $A\longrightarrow A/j(\Ker\rho)$ and
 $A'\longrightarrow A'/j'(\Ker\rho')$ preserve the polarizations up to multiplication
 by some powers of $p$).
 Now we conclude that $y'=(X',\rho')=(X,f^{-1}\circ\rho')=(X,\rho\circ h^{-1})=hy$, as desired.
\end{prf}

Now we can give a proof of our main theorem.

\begin{prf}[of Theorem \ref{thm:main}]\label{proof-main}
 Take $m\ge 0$ as in Proposition \ref{prop:pi-surjection}, and let $\mathcal{S}\subset \Irr(\overline{\mathcal{N}})$ be the subset
 consisting of irreducible components of $\overline{\mathcal{N}}$ which intersect $\overline{\mathcal{N}}{}^{0,m}$.
 Let us observe that $\mathcal{S}$ is a finite set. Take a quasi-compact open subset $U$ of $\overline{\mathcal{N}}$
 containing $\overline{\mathcal{N}}{}^{0,m}$. If $\alpha\in \mathcal{S}$, it intersects $U$ and thus $\alpha\cap U$ is an irreducible
 component of $U$. Moreover the closure of $\alpha\cap U$ in $\overline{\mathcal{N}}$ coincides with $\alpha$.
 Hence there exists an injection $\mathcal{S}\hooklongrightarrow \Irr(U)$. Since $U$ is a scheme of finite type over $\overline{\F}_p$,
 $\Irr(U)$ is a finite set. Thus $\mathcal{S}$ is also a finite set.

 Therefore, it suffices to show that for every $\alpha\in\Irr(\overline{\mathcal{N}})$ there exists $h\in J$ such that
 $h\alpha\in\mathcal{S}$. Fix $y\in \alpha(\overline{\F}_p)$. 
 We can take integers $r$, $m'$ with $m'\ge 0$
 such that $y$ lies in $\overline{\mathcal{N}}{}^{r,m'}(\overline{\F}_p)$. Then $p^{-r}y$ lies in 
 $\overline{\mathcal{N}}{}^{0,m'}(\overline{\F}_p)$. Let $x$ be an arbitrary element of $\Ig_{m'}(\overline{\F}_p)$;
 note that $\Ig_{m'}\neq\varnothing$, as $C$ is non-empty and $\Ig_{m'}\longrightarrow C$ is surjective.
 By Proposition \ref{prop:pi-surjection} ii),
 there exists $(x',y')\in \Ig_m(\overline{\F}_p)\times \overline{\mathcal{N}}{}^{0,m}(\overline{\F}_p)$ such that $\Pi(x,p^{-r}y)=\Pi(x',y')$.
 By Proposition \ref{prop:fiber-orbit}, there exists $h'\in J$ such that $y'=h'p^{-r}y$. Put $h=h'p^{-r}\in J$.
 Then, $h\alpha\in\Irr(\overline{\mathcal{N}})$ intersects $\overline{\mathcal{N}}{}^{0,m}$ at $y'$, and thus $h\alpha\in\mathcal{S}$.
 This completes the proof.
\end{prf}

\section{Applications}\label{sec:applications}
\subsection{Finiteness of $\ell$-adic cohomology of general Rapoport-Zink towers}\label{subsec:general}
Here we continue to use the notation introduced in Section \ref{subsec:statement}.
Let $Z_J$ be the center of $J$.
First we will give a group-theoretic characterization of the quasi-compactness of
$Z_J\backslash\overline{\mathcal{M}}$.

\begin{thm}\label{thm:quasi-compactness}
 Assume that $b$ comes from an abelian variety and $\mathcal{M}$ is non-empty.
 Then, the following are equivalent:
 \begin{enumerate}
  \item[(a)] The group $J$ is compact-mod-center, namely, $J/Z_J$ is compact.
  \item[(b)] The quotient topological space $Z_J\backslash\overline{\mathcal{M}}$ is quasi-compact.
 \end{enumerate}
\end{thm}

\begin{prf}[of Theorem \ref{thm:quasi-compactness} (a)$\implies$(b)]
 Take a compact open subgroup $J^1$ of $J$. Since $J/Z_J$ is compact, the image of $J^1$ in $J/Z_J$
 is a finite index subgroup of $J/Z_J$. Therefore $Z_JJ^1$ is a finite index subgroup of $J$.
 Take a system of representatives $h_1,\ldots,h_k\in J$ of $Z_JJ^1\backslash J$.

 By Theorem \ref{thm:main}, we can choose $\alpha_1,\ldots,\alpha_m\in\Irr(\overline{\M})$
 such that $\Irr(\overline{\M})=\bigcup_{i=1}^mJ\alpha_i=\bigcup_{i=1}^m\bigcup_{j=1}^kZ_JJ^1h_j\alpha_i$.
 By \cite[Proposition 2.3.11]{MR2074714} and the compactness of $J^1$, $J^1h_j\alpha_i$ consists of finitely many elements
 for each $i$ and $j$. Therefore, $Z_J\backslash \overline{\M}$ is covered by the images of finitely many irreducible components
 of $\overline{\M}$. In particular, $Z_J\backslash \overline{\M}$ is quasi-compact.
\end{prf}

To show the converse (b)$\implies$(a), we use the subsequent lemma, which is similar to \cite[Lemma 5.1 iii)]{RZ-LTF}.

\begin{lem}\label{lem:stabilizer-compact}
 Let $T_1$ and $T_2$ be quasi-compact subsets of $\overline{\M}$. Then, the subset
 $\{h\in J\mid hT_1\cap T_2\neq\varnothing\}$ of $J$ is contained in a compact subset of $J$.
\end{lem}

\begin{prf}
 First note that $T_i$ is contained in a finite union of irreducible components of $\overline{\M}$.
 Indeed, take a quasi-compact open subscheme $U$ of $\overline{\M}$ containing $T_i$;
 then, the closure $\overline{\alpha}$ of each $\alpha\in\Irr(U)$ belongs to $\Irr(\overline{\M})$,
 and $T_i$ is contained in $\bigcup_{\alpha\in \Irr(U)}\overline{\alpha}$.
 Therefore, the closure $\overline{T}_i$ of $T_i$ in $\overline{\M}$
 is quasi-compact. 
 Replacing $T_1$, $T_2$ by $\overline{T}_1$, $\overline{T}_2$, we may assume that $T_1$ and $T_2$ are closed.
 For $T=T_1\cup T_2$, $\{h\in J\mid hT_1\cap T_2\neq\varnothing\}$ is contained in $\{h\in J\mid hT\cap T\neq\varnothing\}$.
 Therefore, it suffices to consider the case $T=T_1=T_2$.
 By Lemma \ref{lem:M-N-proper}, we have a $J$-equivariant proper morphism 
 $f\colon \overline{\M}\longrightarrow \overline{\mathcal{N}}$.
 As $\{h\in J\mid hT\cap T\neq\varnothing\}$ is contained in $\{h\in J\mid hf(T)\cap f(T)\neq\varnothing\}$, we may replace
 $\overline{\M}$ by $\overline{\mathcal{N}}$.
 Let $\mathcal{N}_{\GL}$ be the Rapoport-Zink space without any additional structure associated with $\X$, and $J_{\GL}$ the group of
 self-quasi-isogenies of $\X$. Then we have a closed immersion $\mathcal{N}\hooklongrightarrow \mathcal{N}_{\GL}$ compatible with
 the inclusion $J\hooklongrightarrow J_{\GL}$. Hence we may replace $\overline{\mathcal{N}}$ by 
 $\overline{\mathcal{N}}_{\!\GL}=(\mathcal{N}_{\GL})^{\mathrm{red}}$ and $J$ by $J_{\GL}$.
 In this case, the claim is essentially proved in the proof of \cite[Proposition 2.34]{MR1393439}.
\end{prf}

\begin{cor}\label{cor:stabilizer-compact-open}
 \begin{enumerate}
  \item Let $T$ be a non-empty quasi-compact subset of $\overline{\M}$.
	Then, the subset $\{h\in J\mid hT=T\}$ of $J$ is open and compact.
  \item Let $T$ be a finite union of irreducible components of $\overline{\M}$.
	Then, the subset $\{h\in J\mid hT\cap T\neq \varnothing\}$ of $J$ is open and compact.
 \end{enumerate}
\end{cor}

\begin{prf}
 i) It can be proved in the same way as in \cite[Lemma 5.1 iv)]{RZ-LTF}.
 
 We prove ii). The openness follows from i). Write $T=\alpha_1\cup\cdots\cup\alpha_m$, where $\alpha_i\in\Irr(\overline{\M})$.
 Consider the subset $S$ of $\Irr(\overline{\M})$ consisting of $\beta$ satisfying $\beta\cap T\neq \varnothing$.
 It is a finite set (see the proof of Theorem \ref{thm:main} in p.~\pageref{proof-main}).
 Write $S=\{\beta_1,\ldots,\beta_k\}$. Then, we have
 $\{h\in J\mid hT\cap T\neq \varnothing\}=\bigcup_{i=1}^m\bigcup_{j=1}^k\{h\in J\mid h\alpha_i=\beta_j\}$.
 Hence it suffices to show that $\{h\in J\mid h\alpha_i=\beta_j\}$ is compact for every $i$ and $j$.
 We may assume that there exists $h_0\in J$ such that $h_0\alpha_i=\beta_j$. In this case, we have
 $\{h\in J\mid h\alpha_i=\beta_j\}=h_0J_{\alpha_i}$, where $J_{\alpha_i}=\{h\in J\mid h\alpha_i=\alpha_i\}$.
 By i), $J_{\alpha_i}$ is compact, and thus $h_0J_{\alpha_i}$ is also compact. This concludes the proof.
\end{prf}

\begin{prf}[of Theorem \ref{thm:quasi-compactness} (b)$\implies$(a)]
 Assume that $Z_J\backslash \overline{\M}$ is quasi-compact. Then, there exist finitely many irreducible components
 $\alpha_1,\ldots,\alpha_m$ whose images $\overline{\alpha}_1,\ldots,\overline{\alpha}_m$ cover $Z_J\backslash \overline{\M}$.
 Put $T=\alpha_1\cup\cdots\cup\alpha_m$ and $K_T=\{h\in J\mid hT\cap T\neq\varnothing\}$.
 Let us observe that $J=Z_JK_T$. Take any element $h\in J$.
 Since $Z_J\backslash\overline{\M}$ is non-empty,
 we can find $x\in T$. Then there exists $z\in Z_J$ such that $zhx\in T$.
 Thus we have $zh\in K_T$ and $h\in Z_JK_T$, as desired.

 By Corollary \ref{cor:stabilizer-compact-open} ii), $K_T$ is compact. Hence $J/Z_J$ is also compact.
\end{prf}

Next, we will prove a finiteness result on the $\ell$-adic cohomology of the Rapoport-Zink tower.
Before stating the result, we recall some definitions. For a precise description,
see \cite[Chapter 5]{MR1393439}.

We denote the rigid generic fiber $t(\M)_\eta$ of $\M$ by $M$
(note that in this paper all rigid spaces are considered as adic spaces;
\cf \cite{MR1306024}, \cite{MR1734903}).
The universal object on $\M$ induces a system of \'etale $p$-divisible groups $\{\widetilde{X}_L\}_{L\in\mathscr{L}}$
on $M$. 
Assume $M$ is non-empty, and fix a point $x_0$ of $M$ and a geometric point $\overline{x}_0$
lying over $x_0$.
Put $V'=V_p\widetilde{X}_{L,\overline{x}_0}$, which is independent of $L\in\mathscr{L}$. 
Then, the universal quasi-polarization $p_L\colon \widetilde{X}_L\longrightarrow (\widetilde{X}_{{L^\vee}})^\vee$
(\cf the condition (c) in Definition \ref{defn:RZ-space}),
well-defined up to $\Q_p^\times$-multiplication, 
induces an alternating bilinear pairing
$\langle\ ,\ \rangle'\colon V'\times V'\longrightarrow \Q_p$, which is well-defined up to 
$\Q_p^\times$-multiplication (here we choose an isomorphism $\Q_p(1)\cong \Q_p$, but the choice
does not affect the $\Q_p^\times$-orbit of the pairing).
Let $\mathbf{G}'$ be the algebraic group over $\Q_p$ consisting of $B$-linear automorphisms
of $V'$ which preserve $\langle\ ,\ \rangle'$ up to a scalar multiple.
By the comparison theorem for $p$-divisible groups,
there exists a $B\otimes_{\Q_p}B_{\mathrm{crys}}$-linear isomorphism
$V\otimes_{\Q_p}B_{\mathrm{crys}}\yrightarrow{\cong} V'\otimes_{\Q_p}B_{\mathrm{crys}}$
which maps $\langle\ ,\ \rangle$ to a $\Q_p^\times$-multiple of
$\langle\ ,\ \rangle'$. Thus, there exists a $B\otimes_{\Q_p}\overline{\Q}_p$-linear isomorphism
$V\otimes_{\Q_p}\overline{\Q}_p\yrightarrow{\cong} V'\otimes_{\Q_p}\overline{\Q}_p$
which maps $\langle\ ,\ \rangle$ to a $\Q_p^\times$-multiple of
$\langle\ ,\ \rangle'$.
Such a pair $(V',\langle\ ,\ \rangle')$ is classified by $H^1(\Q_p,\mathbf{G})$.
Let $\xi\in H^1(\Q_p,\mathbf{G})$ be the element corresponding to the isomorphism class of 
$(V',\langle\ ,\ \rangle')$.
Then, $\mathbf{G}'$ is the inner form of $\mathbf{G}$ corresponding to the image of $\xi$
under the map $H^1(\Q_p,\mathbf{G})\longrightarrow H^1(\Q_p,\mathbf{G}^{\mathrm{ad}})$.
If $\mathbf{G}$ is connected, there is a formula that describes $\xi$ by means of $b$ and $\mu$;
see \cite[Proposition 1.20]{MR1393439} and \cite[\S 4.5.2]{MR1600395}.
In particular, $\xi$ is independent of the choice of $\overline{x}_0$ in this case.

For $L\in \mathscr{L}$, we denote the $\mathcal{O}_B$-lattice $T_p\widetilde{X}_{L,\overline{x}_0}\subset V'$ by $L'$. Then, $\mathscr{L}'=\{L'\mid L\in\mathscr{L}\}$ is a self-dual multi-chain of
$\mathcal{O}_B$-lattices of $V'$.
We denote by $K'_{\!\mathscr{L}'}$ the subgroup of $G'=\mathbf{G}'(\Q_p)$
consisting of $g$ with $gL'=L'$ for every $L'\in\mathscr{L}'$.
It is compact and open in $G'$.
For an open subgroup $K'$ of $K'_{\!\mathscr{L}'}$,
let $M_{K'}$ be the rigid space over $M$ classifying $K'$-level structures
on $\{\widetilde{X}_L\}_{L\in\mathscr{L}}$. Namely, for a connected rigid space $S$ over $M$,
a morphism $S\longrightarrow M_{K'}$ over $M$ is functorially in bijection with 
a $\pi_1(S,\overline{x})$-invariant $K'$-orbit of systems
of $\mathcal{O}_B$-linear isomorphisms $\{\eta_L\colon L'\yrightarrow{\cong} T_p\widetilde{X}_{L,\overline{x}}\}_{L\in\mathscr{L}}$ such that
\begin{itemize}
 \item $\eta_L$ is compatible with respect to $L$,
 \item and $\eta_L$ maps $\langle\ ,\ \rangle'$ to a $\Q_p^\times$-multiple of
       the alternating bilinear pairing on $V_p\widetilde{X}_{L,\overline{x}}$
       induced from the universal quasi-polarization.
\end{itemize}
Here $\overline{x}$ is a geometric point in $S$; the set of $\pi_1(S,\overline{x})$-invariant
$K'$-orbits of $\{\eta_L\}_{L\in\mathscr{L}}$ is essentially independent of the choice of $\overline{x}$.
We can easily observe that $M_{K'}$ is finite \'etale over $M$. 
If $\mathbf{G}$ is connected, $M_{K'}\longrightarrow M$ is surjective.
The action of $J$ on $M$ naturally lifts to an action on $M_{K'}$.

Varying $K'$, we obtain a projective system $\{M_{K'}\}_{K'\subset K'_{\!\mathscr{L}'}}$
of \'etale coverings over $M$. As usual, we can define a Hecke action of $G'$ on this tower;
see \cite[5.34]{MR1393439}.
By definition, the tower $\{M_{K'}\}_{K'}$ a priori depends on the point $\overline{x}_0$.
However, in fact it is known that $\{M_{K'}\}_{K'}$ depends only on the class
$\xi\in H^1(\Q_p,\mathbf{G})$ (\cf \cite[5.39]{MR1393439}).

Fix a prime number $\ell\neq p$ and consider the compactly supported $\ell$-adic cohomology:
\[
 H^i_c(M_{K'})=H^i_c(M_{K'}\otimes_{\breve{E}}\overline{\breve{E}},\overline{\Q}_\ell),\quad
 H^i_c(M_{\xi,\infty})=\varinjlim_{K'\subset K'_{\!\mathscr{L}'}}H^i_c(M_{K'}).
\]
Then, $H^i_c(M_{\xi,\infty})$ becomes a $G'\times J$-representation.
The action of $G'$ is clearly smooth. By \cite[Corollaire 4.4.7]{MR2074714}, 
the action of $J$ is also smooth. In fact, it is also known that the Weil group $W_E$
naturally acts on $H^i_c(M_{\xi,\infty})$. These three actions are expected to 
be closely related to the local Langlands correspondence (see \cite{MR1403942}).

In the sequel, we will prove the following fundamental finiteness result on the representation
$H^i_c(M_{\xi,\infty})$.

\begin{thm}\label{thm:finiteness}
 Assume that $b$ comes from an abelian variety.
 Then, for every integer $i\ge 0$ and every compact open subgroup $K'$ of $G'$, the $K'$-invariant part $H^i_c(M_{\xi,\infty})^{K'}$
 is finitely generated as a $J$-representation.
\end{thm}

In the unramified case, this theorem is proved in \cite[Proposition 4.4.13]{MR2074714}.
The strategy of our proof is similar. First we will show that $\M$ is locally algebraizable
in the sense of \cite[Definition 3.19]{formalnearby}.
Fix $\widetilde{B}$ and $(A_0,\lambda_0,\iota_0)$ as in Definition \ref{defn:comes-from-AV}.

\begin{lem}
 We can replace $(A_0,\lambda_0,\iota_0)$ so that the following conditions are satisfied:
 \begin{itemize}
  \item there exists an order $\widetilde{\mathcal{O}}$ of $\widetilde{B}$ which is
	contained in $\iota_0^{-1}(\End(A_0))$, stable under $*$
	and satisfies $\widetilde{\mathcal{O}}\otimes_\Z\Z_p=\mathcal{O}_B$.
  \item there exists a finite extension $F$ of $K_0=\Frac W(\overline{\F}_p)$ such that
	$(A_0,\lambda_0,\iota_0)$ lifts to an object $(\widetilde{A}_0,\widetilde{\lambda}_0,\widetilde{\iota}_0)$ over $\mathcal{O}_F$.
 \end{itemize}
\end{lem}

\begin{prf}
 Take a finite extension $F$ of $\breve{E}$ such that $x_0\in M(F)$.
 Then, the point $x_0$ corresponds to an object 
 $\{(\widetilde{X}_L,\widetilde{\iota}_L,\widetilde{\rho}_L)\}_{L\in\mathscr{L}}$ over $\mathcal{O}_F$
 as in Definition \ref{defn:RZ-space}. Let $\{(X_L,\iota_L,\rho_L)\}_{L\in\mathscr{L}}$
 be the element in $\M(\overline{\F}_p)$ that is obtained as the reduction of
 $\{(\widetilde{X}_L,\widetilde{\iota}_L,\widetilde{\rho}_L)\}_{L\in\mathscr{L}}$. 

 Fix $L\in\mathscr{L}$ such that $L\subset L^\vee$.
 There exist an abelian variety $A_0'$ over $\overline{\F}_p$ and a $p$-quasi-isogeny
 $\phi\colon A_0\longrightarrow A_0'$ such that $\phi[p^\infty]\colon A_0[p^\infty]\longrightarrow A'_0[p^\infty]$ 
 can be identified with $\rho_L\colon \X\longrightarrow X_L$.
 By the conditions (a), (c) in Definition \ref{defn:RZ-space}, 
 there exists an integer $m$ such that $(\rho_L^\vee)^{-1}\circ p^m\lambda_0\circ \rho_L^{-1}$
 gives an isogeny $X_L\longrightarrow X_L^\vee$ of height $\log_p\#(L^\vee/L)$. 
 Consider the quasi-isogeny $\lambda_0'=(\phi^{\vee})^{-1}\circ p^m\lambda_0\circ \phi^{-1}\colon A_0'\longrightarrow A_0'^\vee$.
 Passing to $p$-divisible groups, we can easily observe that it is a polarization on $A_0'$.
 On the other hand, let $\iota_0'$ be the composite of
 $\widetilde{B}\yrightarrow{\iota_0} \End(A_0)\otimes_\Z\Q\yrightarrow{(*)}\End(A'_0)\otimes_\Z\Q$,
 where $(*)$ is the isomorphism induced by $\phi$.
 Then, by the same way as in the proof of Lemma \ref{lem:replace-AV}, we can observe that
 $(A'_0,\lambda'_0,\iota'_0)$ satisfies the first condition in the lemma.
 By construction the polarized $B$-isocrystal associated to $(A'_0,\lambda'_0,\iota'_0)$
 is isomorphic to $N_b$.
 Hence we may replace $(A_0,\lambda_0,\iota_0)$ by $(A'_0,\lambda'_0,\iota'_0)$.

 Finally, by the Serre-Tate theorem we obtain a formal lifting 
 $(\widetilde{A}_0,\widetilde{\lambda}_0,\widetilde{\iota}_0)$ to $\mathcal{O}_F$
 corresponding to $(\widetilde{X}_L,\widetilde{\iota}_L)$.
 The existence of the polarization $\widetilde{\lambda}_0$ tells us that $\widetilde{A}_0$ is
 algebraizable.
 This concludes the proof.
\end{prf}
We take $(A_0,\lambda_0,\iota_0)$ and $\widetilde{\mathcal{O}}$ as in the lemma above, 
and fix a lift $(\widetilde{A}_0,\widetilde{\lambda}_0,\widetilde{\iota}_0)$
over $\mathcal{O}_F$.
Let $\mathbf{I}$ be the algebraic group over $\Q$ consisting of $\widetilde{\mathcal{O}}$-linear
self-quasi-isogenies of $A_0$ preserving $\lambda_0$ up to a scalar multiple.
The functor of taking $p$-divisible groups induces an injection
$\mathbf{I}(\Q)\hooklongrightarrow \mathbf{J}(\Q_p)$.

We fix an embedding $\mathcal{O}_F\hooklongrightarrow \C$ and take the base change
$(\widetilde{A}_{0,\C},\widetilde{\lambda}_{0,\C},\widetilde{\iota}_{0,\C})$ of
$(\widetilde{A}_0,\widetilde{\lambda}_0,\widetilde{\iota}_0)$ under this embedding. 
We put $\Lambda=H_1(\widetilde{A}_{0,\C},\Z)$ and $W=H_1(\widetilde{A}_{0,\C},\Q)$.
The $\widetilde{\mathcal{O}}$-action $\widetilde{\iota}_{0,\C}$ on $\widetilde{A}_{0,\C}$
makes $\Lambda$ (resp.\ $W$)
an $\widetilde{\mathcal{O}}$-module (resp.\ $\widetilde{B}$-module).
The polarization $\widetilde{\lambda}_{0,\C}$ induces a $*$-Hermitian alternating bilinear pairing
$\langle\ ,\ \rangle_{\widetilde{\lambda}_{0,\C}}\colon \Lambda\times \Lambda\longrightarrow \Z$.
Clearly we have an isomorphism $\Lambda\otimes_{\Z}\widehat{\Z}^p\cong \varprojlim_{(n,p)=1}A_0[n](\overline{\F}_p)$ compatible with additional structures,
where we write $\widehat{\Z}^p=\varprojlim_{(n,p)=1}\Z/n\Z$.

We denote by $\mathbf{H}$ the algebraic group over $\Q$ consisting of
$\widetilde{B}$-linear automorphisms of $W$
which preserve $\langle\ ,\ \rangle_{\widetilde{\lambda}_{0,\C}}$ up to a scalar multiple.
We have a natural injection $\mathbf{I}(\Q)\hooklongrightarrow \mathbf{H}(\A_f^p)$, where
$\A_f^p=\widehat{\Z}^p\otimes_\Z\Q$.

\begin{rem}
 If $\mathbf{G}$ is connected, we can also compare $V\otimes_{\Q_p}K_0=\mathbb{D}(A_0[p^\infty])_\Q$
 with $W\otimes_\Q K_0$.
 By the comparison result between de Rham and crystalline theories (\cf \cite{MR700767}),
 we have a canonical isomorphism $V\otimes_{\Q_p}\C\cong H_1(\widetilde{A}_{0,\C},\Q)\otimes_\Q\C$.
 Therefore, by Steinberg's theorem $H^1(K_0,\mathbf{G})=1$, we have an isomorphism
 $V\otimes_{\Q_p}K_0\cong H_1(\widetilde{A}_{0,\C},\Q)\otimes_\Q K_0$ 
 compatible with various additional structures. However, we do not need this result.
\end{rem}

We will use the following moduli space, which is a slight variant of \cite[Proposition 6.9]{MR1393439}.

\begin{defn}
 For a compact open subgroup $K^p$ of $\mathbf{H}(\A_f^p)$,
 let $\mathfrak{Y}_{K^p}$ be the functor from the category of
 $\mathcal{O}_{\breve{E}}$-schemes to $\mathbf{Set}$
 that associates $S$ to the set of isomorphism classes of
 $(A,\overline{\lambda},\alpha)$ where
 \begin{itemize}
 \item $A=\{A_L\}$ is an $\mathscr{L}$-set of $\widetilde{\mathcal{O}}$-abelian schemes over $S$
       (\cf \cite[Definition 6.5]{MR1393439}; recall that we are working on the category of
       $\widetilde{\mathcal{O}}$-abelian schemes up to isogeny of order prime to $p$)
	such that for each $L\in\mathscr{L}$ the determinant condition
	$\det\nolimits_{\mathcal{O}_S}(a;\Lie A_L)=\det\nolimits_K(a;V_0)$
	in $a\in\widetilde{\mathcal{O}}$ holds,
 \item $\overline{\lambda}$ is a $\Q$-homogeneous principal polarization of $A$
       (\cf \cite[Definition 6.7]{MR1393439}),
 \item and $\alpha$ is a $K^p$-level structure
       \[
       \alpha\colon H_1(A,\A^p_f)\yrightarrow{\cong} W\otimes_\Q\A_f^p\pmod{K^p}
       \]
       which is $\widetilde{B}\otimes_\Q\A_f^p$-linear and preserves the pairings up to $(\A_f^p)^\times$-multiplication (for the definition of $H_1(A,\A^p_f)$, see \cite[6.8]{MR1393439}).
 \end{itemize}
 As in \cite[p.~279]{MR1393439}, we can easily show that the functor $\mathfrak{Y}_{K^p}$ is
 represented by a quasi-projective scheme over $\mathcal{O}_{\breve{E}}$.
\end{defn}

The difference between \cite[Definition 6.9]{MR1393439} and this definition is that
in our case $\mathscr{L}$ is a self-dual multi-chain of $\mathcal{O}_B$-lattices in $V$,
not in $W\otimes_\Q\Q_p$.
However, the proof of \cite[Theorem 6.23]{MR1393439} can be applied to our case without any change,
and we obtain the following $p$-adic uniformization result.

\begin{prop}
 \begin{enumerate}
  \item Let $\widehat{\mathfrak{Y}}_{K^p}$ be the $p$-adic completion of $\mathfrak{Y}_{K^p}$.
	We have a morphism of pro-formal schemes over $\mathcal{O}_{\breve{E}}$:
	\[
	 \Theta\colon \M\times \mathbf{H}(\A_f^p)/K^p\longrightarrow \widehat{\mathfrak{Y}}_{K^p}.
	\]
  \item The group $\mathbf{I}(\Q)$ is discrete in $\mathbf{J}(\Q_p)\times \mathbf{H}(\A_f^p)$.
  \item If $K^p$ is small enough, the quotient $\mathbf{I}(\Q)\backslash\M\times \mathbf{H}(\A_f^p)/K^p$
	is a formal scheme. It is a countable disjoint union of spaces of the form $\Gamma\backslash\M$,
	where $\Gamma$ is a subgroup of $\mathbf{J}(\Q_p)$ of the form
	$(\mathbf{J}(\Q_p)\times hK^ph^{-1})\cap \mathbf{I}(\Q)$ with $h\in \mathbf{H}(\A_f^p)$.
	Such a subgroup $\Gamma$ is torsion-free and discrete in $\mathbf{J}(\Q_p)$.
  \item Let $\mathcal{T}$ be the set of closed subsets of $\mathfrak{Y}_{K^p}\otimes_{\mathcal{O}_{\breve{E}}}\overline{\F}_p$ which are the images of irreducible components of $\M\times \mathbf{H}(\A_f^p)/K^p$
	under $\Theta$. Then, each $T\in \mathcal{T}$ intersects only finitely many members of $\mathcal{T}$.
  \item The morphism $\Theta$ induces an isomorphism of formal schemes over $\mathcal{O}_{\breve{E}}$:
	\[
	 \Theta\colon \mathbf{I}(\Q)\backslash\M\times \mathbf{H}(\A_f^p)/K^p\yrightarrow{\cong} (\mathfrak{Y}_{K^p})_{/\mathcal{T}}.
	\]
	The target $(\mathfrak{Y}_{K^p})_{/\mathcal{T}}$ is the formal completion of
	$\mathfrak{Y}_{K^p}$
	along $\mathcal{T}$; see \cite[6.22]{MR1393439} for a precise definition.
 \end{enumerate}
\end{prop}
The claim i) is proved in \cite[Theorem 6.21]{MR1393439} (under the setting therein).
The remaining statements are found in \cite[Theorem 6.23]{MR1393439} and its proof.

By the same argument as in the proof of \cite[Corollaire 3.1.4]{MR2074714}, we obtain the following
algebraization result.

\begin{cor}\label{cor:RZ-loc-alg}
 For every quasi-compact open formal subscheme $\mathcal{U}$ of $\M$, we can find a compact open subgroup
 $K^p$ of $\mathbf{H}(\A_f^p)$ and a locally closed subscheme $Z$ of
 $\mathfrak{Y}_{K^p}\otimes_{\mathcal{O}_{\breve{E}}}\overline{\F}_p$ such that
 $\mathcal{U}$ is isomorphic to the formal completion of $\mathfrak{Y}_{K^p}$ along $Z$.
 In particular, $\M$ is locally algebraizable in the sense of \cite[Definition 3.19]{formalnearby}.
\end{cor}

Next we consider level structures at $p$. We write $\mathfrak{Y}_{K^p,\eta}$ for the generic fiber
$\mathfrak{Y}_{K^p}\otimes_{\mathcal{O}_{\breve{E}}}\breve{E}$, and $\widetilde{A}=\{\widetilde{A}_L\}$
the universal $\mathscr{L}$-set of $\widetilde{\mathcal{O}}$-abelian schemes over $\mathfrak{Y}_{K^p,\eta}$.
For an open subgroup $K'\subset K'_{\!\mathscr{L}'}$, consider a $K'$-level structure
\[
 \{L'\yrightarrow{\cong} T_p\widetilde{A}_L \ (\bmod{K'})\}_{L\in\mathscr{L}}
\]
(the definition is similar as that in the definition of $M_{K'}$).
Let $\mathfrak{Y}_{K'K^p,\eta}$ be the scheme classifying such level structures,
which is finite \'etale over $\mathfrak{Y}_{K^p,\eta}$.

The following corollary follows directly from the proof of \cite[Corollaire 3.1.4]{MR2074714}:

\begin{cor}\label{cor:RZ-alg-with-level}
 Let $\mathcal{U}$ be a quasi-compact open formal subscheme of $\M$, and $K'$ be
 an open subgroup of $K'_{\!\mathscr{L}'}$.
 Put $U=t(\mathcal{U})_\eta$ and denote by $U_{K'}$ the inverse image of $U$ under $M_{K'}\longrightarrow M$.
 The locally closed subscheme $Z$ of $\mathfrak{Y}_{K^p}\otimes_{\mathcal{O}_{\breve{E}}}\overline{\F}_p$
 in Corollary \ref{cor:RZ-loc-alg} can be taken so that we have the following cartesian diagram:
 \[
  \xymatrix{%
 U_{K'}\ar[rr]\ar[d]& & (\mathfrak{Y}_{K'K^p,\eta})^{\mathrm{ad}}\ar[d]\\
 U\ar[r]^-{\cong}& t((\mathfrak{Y}_{K^p})_{/Z})_\eta\ar[r]^-{\subset}_-{\mathrm{open}}&(\mathfrak{Y}_{K^p,\eta})^{\mathrm{ad}}\lefteqn{.}
 }
 \]
 Here $(-)^{\mathrm{ad}}$ denotes the associated adic space, and $(-)_{/Z}$ the formal completion along $Z$.
\end{cor}

\begin{cor}\label{cor:finiteness-cohomology}
 Let the setting as in Corollary \ref{cor:RZ-alg-with-level}. Then,
 $H^i_c(U_{K'}\otimes_{\breve{E}}\overline{\breve{E}},\overline{\Q}_\ell)$ is a finite-dimensional
 $\overline{\Q}_\ell$-vector space.
\end{cor}

\begin{prf}
 By Zariski's main theorem, there exists a scheme $\mathfrak{Y}_{K'K^p}$ finite over $\mathfrak{Y}_{K^p}$ that contains $\mathfrak{Y}_{K'K^p,\eta}$ as a dense open subscheme. 
 Since $\mathfrak{Y}_{K'K^p,\eta}$ is finite over $\mathfrak{Y}_{K^p,\eta}$, we conclude that
 the generic fiber of $\mathfrak{Y}_{K'K^p}$ coincides with $\mathfrak{Y}_{K'K^p,\eta}$.
 We write $Z'$ for the inverse image of $Z$ under $\mathfrak{Y}_{K'K^p}\longrightarrow \mathfrak{Y}_{K^p}$. Then, by Corollary \ref{cor:RZ-alg-with-level}, we have an isomorphism
 $U_{K'}\cong t((\mathfrak{Y}_{K'K^p})_{/Z'})_\eta$.
 Thus the finiteness follows from \cite[Lemma 3.13 i)]{MR1620118}, \cite[Corollary 2.11]{MR1620114}
 and \cite[Lemma 3.3 (ii)]{MR1626021}.
\end{prf}

Now we can give a proof of Theorem \ref{thm:finiteness}.

\begin{prf}[of Theorem \ref{thm:finiteness}]
 The proof is similar as that of \cite[Proposition 4.4.13]{MR2074714}
 (see also \cite[\S 5]{RZ-LTF}).
 We include it for the sake of completeness.

 Recall that every $J$-submodule of a finitely generated smooth $J$-representation is again
 finitely generated. If $\mathbf{J}$ is connected, this result is found in \cite[Remarque 3.12]{MR771671}; the general case is immediately reduced to the connected case.
 In particular we may assume that $K'\subset K'_{\mathscr{L}}$. 
 In this case, we have $H^i_c(M_{\xi,\infty})^{K'}=H^i_c(M_{K'})$. We will show that
 it is finitely generated as a $J$-representation.

 For simplicity, write $\mathcal{I}=\Irr(\overline{\M})$. For $\alpha\in \mathcal{I}$, put
 $\overline{\mathcal{U}}_\alpha=\overline{\M}\setminus\bigcup_{\beta\in \mathcal{I},\alpha\cap\beta=\varnothing}\beta$. It is a quasi-compact open subscheme of $\overline{\M}$.
 Note that there exists only finitely many $\beta\in \mathcal{I}$ with 
 $\overline{\mathcal{U}}_\alpha\cap \overline{\mathcal{U}}_\beta\neq \varnothing$
 (\cf \cite[Corollary 5.2]{RZ-LTF}).
 Let $\mathcal{U}_\alpha$ be
 the open formal subscheme of $\M$ satisfying $(\mathcal{U}_\alpha)^{\mathrm{red}}=\overline{\mathcal{U}}_\alpha$, and $U_\alpha$ the rigid generic fiber of $\mathcal{U}_\alpha$.
 The inverse image $U_{\alpha,K'}$ of $U_\alpha$ under $M_{K'}\longrightarrow M$
 gives an open covering $\{U_{\alpha,K'}\}_{\alpha\in \mathcal{I}}$ of $M_{K'}$.
 For $h\in J$, we have $hU_{\alpha,K'}=U_{h\alpha,K'}$.

 For a finite non-empty subset $\underline{\alpha}=\{\alpha_1,\ldots,\alpha_m\}$ of $\mathcal{I}$,
 we put $U_{\underline{\alpha},K'}=\bigcap_{j=1}^mU_{\alpha_j,K'}$.
 By Corollary \ref{cor:finiteness-cohomology} for $\bigcap_{j=1}^m\mathcal{U}_{\alpha_j}$, 
 $H^i_c(U_{\underline{\alpha},K'})=H^i_c(U_{\underline{\alpha},K'}\otimes_{\breve{E}}\overline{\breve{E}},\overline{\Q}_\ell)$ is a finite-dimensional $\overline{\Q}_\ell$-vector space.
 Put $J_{\underline{\alpha}}=\{h\in J\mid h\underline{\alpha}=\underline{\alpha}\}$.
 By Corollary \ref{cor:stabilizer-compact-open} i), $J_{\underline{\alpha}}$ is
 a compact open subgroup of $J$. It acts smoothly on $H^i_c(U_{\underline{\alpha},K'})$.
 
 For an integer $s\ge 1$, let $\mathcal{I}_s$ be the set of subsets $\underline{\alpha}\subset\mathcal{I}$ such that $\#\underline{\alpha}=s$ and $U_{\underline{\alpha},K'}\neq \varnothing$.
 Note that there exists an integer $N$ such that $\mathcal{I}_s=\varnothing$ for $s>N$.
 Indeed, by Theorem \ref{thm:main}, we can take a finite system of representatives
 $\alpha_1,\ldots,\alpha_k$ of $J\backslash \Irr(\overline{\M})$; then we may take 
 $N$ as the maximum of $\#\{\beta\in\Irr(\overline{\M})\mid \overline{\mathcal{U}}_{\alpha_j}\cap \overline{\mathcal{U}}_\beta\neq\varnothing\}$
 for $1\le j\le k$.

 The group $J$ naturally acts on $\mathcal{I}_s$.
 We will observe that this action has finite orbits.
 Let $\mathcal{I}_s^\sim$ be the subset of $\mathcal{I}^s$ consisting of $(\alpha_1,\ldots,\alpha_s)$
 such that $\alpha_1,\ldots,\alpha_s$ are mutually disjoint and $U_{\{\alpha_1,\ldots,\alpha_s\},K'}\neq\varnothing$. As we have a $J$-equivariant surjection $\mathcal{I}_s^\sim\longrightarrow \mathcal{I}_s$;
 $(\alpha_1,\ldots,\alpha_s)\longmapsto \{\alpha_1,\ldots,\alpha_s\}$, it suffices to show that
 the action of $J$ on $\mathcal{I}_s^\sim$ has finite orbits.
 Consider the first projection $\mathcal{I}_s^\sim\longrightarrow \mathcal{I}$,
 which is obviously $J$-equivariant. The fiber of this map is finite, since
 for each $\alpha\in \mathcal{I}$ there exist only finitely many $\beta\in \mathcal{I}$
 such that $U_{\alpha,K'}\cap U_{\beta,K'}\neq \varnothing$. On the other hand,
 by Theorem \ref{thm:main}, the action of $J$ on $\mathcal{I}$ has finite orbits.
 Hence the action of $J$ on $\mathcal{I}_s^\sim$ also has finite orbits, as desired.
 For each $s\ge 1$, we fix a system of representatives 
 $\underline{\alpha}_{s,1},\ldots,\underline{\alpha}_{s,k_s}$ of $J\backslash \mathcal{I}_s$.

 Consider the \v{C}ech spectral sequence 
 \[
  E_1^{-s,t}=\bigoplus_{\underline{\alpha}\in\mathcal{I}_{s+1}}H^t_c(U_{\underline{\alpha},K'})\Longrightarrow H^{-s+t}_c(M_{K'})
 \]
 with respect to the covering $\{U_{\alpha,K'}\}_{\alpha\in\mathcal{I}}$. This spectral sequence
 is $J$-equivariant, and we can easily see that $E_1^{-s,t}\cong \bigoplus_{j=1}^{k_s}\cInd_{J_{\underline{\alpha}_{s,j}}}^JH^t_c(U_{\underline{\alpha}_{s,j},K'})$ as $J$-representations.
 The right hand side is obviously finitely generated. Thus so are the $E_1$-terms.
 Hence, by the property of finitely generated smooth $J$-representations recalled
 at the beginning of the proof,
 we may conclude that $H^i_c(M_{K'})$ is finitely generated as a $J$-representation for every $i$.
\end{prf}

\subsection{Stronger finiteness of $\ell$-adic cohomology of the basic Rapoport-Zink tower for $\GSp_{2n}$}\label{subsec:GSp}
As in Section \ref{subsec:GU(n,D)}, let $D=\Q_{p^2}[\Pi]$ be a quaternion division algebra over $\Q_p$ and
$*$ an involution on $D$ defined by $\Pi^*=\Pi$ and $a^*=\sigma(a)$ for every $a\in \Q_{p^2}$.
We write $\mathcal{O}_D$ for the maximal order of $D$.

Let us fix an integer $n\ge 1$. We consider non-degenerate alternating bilinear pairings
\begin{itemize}
 \item $\langle\ ,\ \rangle\colon \Q_p^{2n}\times \Q_p^{2n}\longrightarrow \Q_p$; 
       $\langle (x_i),(y_i)\rangle=\sum_{i=1}^n(x_{2i-1}y_{2i}-x_{2i}y_{2i-1})$, and
 \item $\langle\ ,\ \rangle_D\colon D^n\times D^n\longrightarrow \Q_p$;
       $\langle (d_i),(d'_i)\rangle_D=\sum_{i=1}^n\Trd(\varepsilon d_i^*d'_i)$,
\end{itemize}
where $\varepsilon$ is an element of $\Z^\times_{p^2}$ satisfying $\sigma(\varepsilon)=-\varepsilon$.
Note that the lattice $\Z_p^{2n}$ is self-dual with respect to $\langle\ ,\ \rangle$, and
the dual with respect to $\langle\ ,\ \rangle_D$ of the lattice $\mathcal{O}_D^n\subset D^n$ is
equal to $\Pi^{-1}\mathcal{O}_D^n$.
Let $\mathrm{GSp}_{2n}$ denote the symplectic similitude group attached to $(\Q_p^{2n},\langle\ ,\ \rangle)$,
and $\GU(n,D)$ denote the quaternion unitary similitude group attached to $(D^n,\langle\ ,\ \rangle_D)$.
We fix an isomorphism between $(\Q_p^{2n}\otimes_{\Q_p}\Q_{p^2},\langle\ ,\ \rangle)$
and the object attached to $(D^n\otimes_{\Q_p}\Q_{p^2},\langle\ ,\ \rangle_D)$
by the equivalence in Lemma \ref{lem:D-decomposition} i). It gives rise to an inner twist 
$\psi\colon \GU(n,D)\otimes_{\Q_p}\Q_{p^2}\yrightarrow{\cong}\mathrm{GSp}_{2n}\otimes_{\Q_p}\Q_{p^2}$
(see Remark \ref{rem:D-decomp} ii)).

We consider the Rapoport-Zink datum
\[
 \mathcal{D}=(\Q_p,\Z_p,\id,\Q_p^{2n},\langle\ ,\ \rangle,\{p^m\Z_p^{2n}\}_{m\in \Z},b,\mu),
\]
where
\begin{itemize}
 \item $b$ is an element of $\GSp_{2n}(K_0)$ such that the polarized $\Q_p$-isocrystal
       $N_b$ is isomorphic to $\mathbb{D}(E_0[p^\infty]^{\oplus n})_\Q$, where $E_0$ is a supersingular
       elliptic curve over $\overline{\F}_p$ endowed with the standard polarization,	       
 \item and $\mu\colon \mathbb{G}_m\longrightarrow \GSp_{2n}$ is the cocharacter given by
       $\mu(z)=\mathrm{diag}(z,1,\ldots,z,1)$.
\end{itemize}
Under the equivalence in Lemma \ref{lem:D-decomposition} ii), the $p$-polarized $\Q_p$-isocrystal $N_b$
corresponds to the $p^2$-polarized $D$-isocrystal $N_p$ attached to a central element 
$p\in K_0^\times\subset \GU(n,D)(K_0)$
(we can check it by using Remark \ref{rem:GL-Newton} and Lemma \ref{lem:D-isoc-Newton-inj}).
Hence the algebraic group $\mathbf{J}$ in Section \ref{subsec:statement} is identified with $\GU(n,D)$.
We also write $\GU(n,D)$ for the group $\mathbf{J}(\Q_p)$ of $\Q_p$-valued points.

We denote by $\mathcal{M}$ the Rapoport-Zink space attached to $\mathcal{D}$.
As in the previous subsection, we can construct the Rapoport-Zink tower $\{M_K\}_K$,
where $K$ runs through enough small compact open subgroups of $\GSp_{2n}(\Q_p)$.
We call it the basic Rapoport-Zink tower for $\GSp_{2n}$.
Note that $H^1(\Q_p,\GSp_{2n})=1$ implies that $\xi=1$ and $\mathbf{G}'=\GSp_{2n}$.
Therefore, the compactly supported $\ell$-adic cohomology
\[
 H^i_c(M_\infty)=\varinjlim_K H^i_c(M_K\otimes_{\breve{\Q}_p}\overline{\breve{\Q}}_p,\overline{\Q}_\ell)
\]
is endowed with a smooth action of $\GSp_{2n}(\Q_p)\times \GU(n,D)$.

We also introduce another Rapoport-Zink tower, which is dual to $\{M_K\}_K$.
Consider the Rapoport-Zink datum
\[
 \check{\mathcal{D}}=(D,\mathcal{O}_D,*,D^n,\langle\ ,\ \rangle_D,\{\Pi^m\mathcal{O}_D^n\}_{m\in \Z},\check{b},\check{\mu}),
\]
where
\begin{itemize}
 \item $\check{b}$ is an element of $\GU(n,D)(K_0)$ such that the polarized $D$-isocrystal
       $N_{\check{b}}$ is isomorphic to $N_{0,1}^{\oplus n}$, where $N_{0,1}$ is the simple polarized $D$-isocrystal
       of type $(0,1)$ introduced in Definition \ref{defn:simple-D-isoc},
 \item and $\check\mu\colon \mathbb{G}_m\longrightarrow \GU(n,D)\otimes_{\Q_p}K_0$ is the cocharacter given by
       the composite of
       \[
	\mathbb{G}_m\yrightarrow{\mu} \GSp_{2n}\otimes_{\Q_p}K_0\yrightarrow[\cong]{\psi^{-1}} \GU(n,D)\otimes_{\Q_p}K_0.
       \]
\end{itemize}
By Lemma \ref{lem:D-decomposition} ii), the algebraic group $\mathbf{J}$ attached to this Rapoport-Zink datum
is identified with $\GSp_{2n}$.

We denote by $\check{\mathcal{M}}$ the Rapoport-Zink space attached to $\check{\mathcal{D}}$,
and by $\{\check{M}_H\}_H$ the Rapoport-Zink tower lying over $t(\check{\mathcal{M}})_\eta$.
Here $H$ runs through enough small compact open subgroups of $\GU(n,D)$.
Note that $H^1(\Q_p,\GU(n,D))=1$ implies that $\xi=1$ and $\mathbf{G}'=\GU(n,D)$ in this case.
Therefore, the compactly supported $\ell$-adic cohomology
\[
 H^i_c(\check{M}_\infty)=\varinjlim_H H^i_c(\check{M}_H\otimes_{\breve{\Q}_p}\overline{\breve{\Q}}_p,\overline{\Q}_\ell)
\]
is endowed with a smooth action of $\GSp_{2n}(\Q_p)\times \GU(n,D)$.

The following theorem is a consequence of the duality isomorphism proved in
\cite{2017arXiv170906651K}, \cite{2017arXiv171006935C} and \cite{Scholze-Berkeley}:

\begin{thm}\label{thm:duality-isom}
 We have a $\GSp_{2n}(\Q_p)\times \GU(n,D)$-equivariant isomorphism 
 \[
  H^i_c(M_\infty)\cong H^i_c(\check{M}_\infty).
 \]
\end{thm}

\begin{prf}
 For $b'\in \GSp_{2n}(K_0)$ and a cocharacter $\mu'\colon \mathbb{G}_m\longrightarrow \GSp_{2n,\overline{K}_0}$,
 one can construct a locally spatial diamond $M(b',\mu')$ over $\C_p$ 
 (see \cite[Definition 4.8.3]{2017arXiv170906651K},
 \cite[Theorem 3.3]{2017arXiv171006935C}, \cite[\S 23.1]{Scholze-Berkeley}).
 We briefly recall the construction of it.
 For a perfectoid $\C_p$-algebra $R$ with tilt $R^\flat$, the Fargues-Fontaine curve $X_{R^\flat}$ and
 a Cartier divisor $D_R\hooklongrightarrow X_{R^\flat}$ are naturally attached.
 Each element $b''\in \GSp_{2n}(K_0)$ determines a $\GSp_{2n}$-bundle $\mathcal{E}_{b''}$ on $X_{R^\flat}$.
 The diamond $M(b',\mu')$ is the sheafification of the functor that attaches a perfectoid affinoid $\Spa(R,R^+)$ over $\C_p$ to
 the set of isomorphism classes of modifications
 $\mathcal{E}_{b'}\vert_{X_{R^\flat}\setminus D_R}\yrightarrow{\cong} \mathcal{E}_1\vert_{X_{R^\flat}\setminus D_R}$
 which are fiberwise bounded by $\mu'$.
 It is equipped with an action of $\GSp_{2n}(\Q_p)\times \mathbf{J}_{b'}(\Q_p)$,
 where $\mathbf{J}_{b'}$ denotes the algebraic group of automorphisms of the polarized $\Q_p$-isocrystal $N_{b'}$.
 If $(b',\mu')=(b,\mu)$, $M(b,\mu)$ is known to be a perfectoid space, which is characterized by
 the property $M(b,\mu)\sim\varprojlim_K M_K\otimes_{\breve{\Q}_p}\C_p$
 (see \cite[Corollary 24.3.5]{Scholze-Berkeley}).
 Let $\phi$ be the cocharacter $\mathbb{G}_m\longrightarrow\GSp_{2n}$; $z\longmapsto \mathrm{diag}(z,\ldots,z)$
 and consider the case $(b',\mu')=(p^{-1}b,\phi^{-1}\mu)$.
 By definition, it is easy to observe that there exists a modification
 $\mathcal{E}_b\vert_{X_{R^\flat}\setminus D_R}\yrightarrow{\cong} \mathcal{E}_{p^{-1}b}\vert_{X_{R^\flat}\setminus D_R}$
 with relative position $\phi$. This gives an isomorphism $M(b,\mu)\cong M(p^{-1}b,\phi^{-1}\mu)$,
 which is $\GSp_{2n}(\Q_p)\times \GU(n,D)$-equivariant under the identification $\GU(n,D)\cong \mathbf{J}_b=\mathbf{J}_{p^{-1}b}$.
 Since $\phi^{-1}\mu$ is conjugate to $\mu^{-1}$, we obtain an isomorphism $M(b,\mu)\cong M(p^{-1}b,\mu^{-1})$.
 By similar construction, we also have a perfectoid space $M(\check{b},\check{\mu})$
 equipped with an action of $\GSp_{2n}(\Q_p)\times \GU(n,D)$, which is characterized by 
 the property $M(\check{b},\check{\mu})\sim\varprojlim_H \check{M}_H\otimes_{\breve{\Q}_p}\C_p$.

 As in the proof of Lemma \ref{lem:D-isoc-Newton-inj}, for an algebraic group $\mathbf{G}$ over $\Q_p$,
 we write $\mathbf{B}(\mathbf{G})$ for the set of $\sigma$-conjugacy classes in $\mathbf{G}(K_0)$.
 For $g\in \mathbf{G}(K_0)$, we write $[g]$ for the $\sigma$-conjugacy class of $g$.
 By Lemma \ref{lem:D-decomposition} ii) and Remark \ref{rem:D-decomp} i), we have a bijection
 $\mathbf{B}(\GU(n,D))\yrightarrow{\cong}\mathbf{B}(\GSp_{2n})$ characterized by the following property:
 an element $[b_1]\in \mathbf{B}(\GU(n,D))$ is mapped to $[b_0]\in \mathbf{B}(\GSp_{2n})$ if
 the polarized $D$-isocrystal $N_{b_1}$ corresponds to the polarized $\Q_p$-isocrystal $N_{b_0}$ under the equivalence
 in Lemma \ref{lem:D-decomposition} ii). We can easily check that it coincides with the bijection in 
 \cite[\S 5.1]{2017arXiv171006935C} for $(G,b)=(\GU(n,D),\check{b})$.
 Further, by using Remark \ref{rem:GL-Newton} and Lemma \ref{lem:D-isoc-Newton-inj}, 
 we can observe that the image of $[1]\in \mathbf{B}(\GU(n,D))$
 is equal to $[p^{-1}b]\in \mathbf{B}(\GSp_{2n})$.
 Therefore, the duality isomorphism (see \cite[Proposition 4.9.1]{2017arXiv170906651K}, \cite[\S 5.1]{2017arXiv171006935C} and \cite[Corollary 23.2.3]{Scholze-Berkeley})
 tells us that there exists a $\GSp_{2n}(\Q_p)\times \GU(n,D)$-equivariant isomorphism
 $M(p^{-1}b,\mu^{-1})\cong M(\check{b},\check{\mu})$. Hence we have $M(b,\mu)\cong M(\check{b},\check{\mu})$.
 Together with \cite[Proposition 5.4 (iii)]{Scholze-perfectoid-survey}, we obtain
 \[
  H^i_c(M_\infty)\cong H^i_c(M(b,\mu),\overline{\Q}_\ell)^{\mathrm{sm}}\cong H^i_c(M(\check{b},\check{\mu}),\overline{\Q}_\ell)^{\mathrm{sm}}\cong H^i_c(\check{M}_\infty),
 \]
 where $(-)^{\mathrm{sm}}$ denotes the $\GSp_{2n}(\Q_p)\times \GU(n,D)$-smooth part.
 This concludes the proof.
\end{prf}

Together with Theorem \ref{thm:finiteness}, we obtain the following corollary.

\begin{cor}\label{cor:finiteness+duality}
 For every $i\ge 0$, the $\GSp_{2n}(\Q_p)\times \GU(n,D)$-representation $H^i_c(M_\infty)$ satisfies
 the following:
 \begin{enumerate}
  \item For every compact open subgroup $K$ of $\GSp_{2n}(\Q_p)$,
	$H^i_c(M_\infty)^K$ is a finitely generated $\GU(n,D)$-representation.
  \item For every compact open subgroup $H$ of $\GU(n,D)$,
	$H^i_c(M_\infty)^H$ is a finitely generated $\GSp_{2n}(\Q_p)$-representation.
 \end{enumerate}
\end{cor}

\begin{prf}
 Note that $b$ comes from an abelian variety by definition,
 and $\check{b}$ comes from an abelian variety by Proposition \ref{prop:GU(n,D)-comes-from-AV}.
 Therefore i) is a direct consequence of Theorem \ref{thm:finiteness},
 and ii) follows from Theorem \ref{thm:finiteness} and Theorem \ref{thm:duality-isom}.
\end{prf}

We denote by $\mathcal{D}_c$ the convolution algebra of distributions on $\GSp_{2n}(\Q_p)$ with
compact support. For a compact open subgroup $K$ of $\GSp_{2n}(\Q_p)$, the idempotent $e_K$ of
$\mathcal{D}_c$ is naturally attached. For a (left or right) $\mathcal{D}_c$-module $V$, we put
$V^{\text{$\mathcal{D}_c$-sm}}=\varinjlim_K e_KV$, where $K$ runs through compact open subgroups
of $\GSp_{2n}(\Q_p)$. It is a smooth representation of $\GSp_{2n}(\Q_p)$.

Let $i,r\ge 0$ be integers. For an irreducible smooth representation $\rho$ of $\GU(n,D)$,
we consider the $r$th Ext group $\Ext^r_{\GU(n,D)}(H^i_c(M_\infty),\rho)$ in the category of
smooth $\GU(n,D)$-representations. Since a smooth representation of $\GSp_{2n}(\Q_p)$ can be
regarded as a left $\mathcal{D}_c$-module, $\Ext^r_{\GU(n,D)}(H^i_c(M_\infty),\rho)$ carries
a right $\mathcal{D}_c$-module structure. Hence we can consider 
$\Ext^r_{\GU(n,D)}(H^i_c(M_\infty),\rho)^{\text{$\mathcal{D}_c$-sm}}$, which is a smooth representation
of $\GSp_{2n}(\Q_p)$. If $r=0$, we can verify that 
$\Hom_{\GU(n,D)}(H^i_c(M_\infty),\rho)^{\text{$\mathcal{D}_c$-sm}}=\Hom_{\GU(n,D)}(H^i_c(M_\infty),\rho)^{\text{sm}}$, where $(-)^{\mathrm{sm}}$ denotes the set of $\GSp_{2n}(\Q_p)$-smooth vectors.

The main result of this subsection is the following:

\begin{thm}\label{thm:coh-fin-length}
 For every integers $i,r\ge 0$ and every irreducible smooth representation $\rho$ of $\GU(n,D)$,
 the $\GSp_{2n}(\Q_p)$-representation
 $\Ext^r_{\GU(n,D)}(H^i_c(M_\infty),\rho)^{\text{$\mathcal{D}_c$-sm}}$ has finite length.
 In particular, the $\GSp_{2n}(\Q_p)$-representation
 $\Hom_{\GU(n,D)}(H^i_c(M_\infty),\rho)^{\mathrm{sm}}$ has finite length.
\end{thm}

We deduce Theorem \ref{thm:coh-fin-length} from Corollary \ref{cor:finiteness+duality} by using
the theory of the Bernstein decomposition (\cf \cite{MR771671}, \cite[Chapitre VI]{MR2567785}).
To introduce some notation, let us consider a general connected reductive group $\mathbf{G}$ over $\Q_p$.
We put $G=\mathbf{G}(\Q_p)$.
Let $\mathrm{Rep}(G)$ be the category of smooth representations of $G$, $\Omega(G)$ the set of equivalence classes of cuspidal data
(\cf \cite[Definition VI.5.2]{MR2567785}), and $\mathcal{B}(G)$ the set of connected components of $\Omega(G)$.
For an element $\mathfrak{s}\in \mathcal{B}(G)$, the subcategory $\mathrm{Rep}(G)_{\mathfrak{s}}$ of $\mathrm{Rep}(G)$ is naturally attached
and the category $\mathrm{Rep}(G)$ is decomposed into the direct product
$\prod_{\mathfrak{s}\in\mathcal{B}(G)}\mathrm{Rep}(G)_{\mathfrak{s}}$.

Here we say that a compact open subgroup $K$ of $G$ is nice
if it is a normal subgroup of a special maximal compact subgroup
and has the Iwahori decomposition.
Such subgroups form a system of fundamental neighborhoods of the unit in $G$
(see \cite[Th\'eor\`eme V.5.2]{MR2567785}).
For a nice compact open subgroup $K$ of $G$, we denote by $\mathrm{Rep}(G)_K$ the subcategory of
$\mathrm{Rep}(G)$ consisting of smooth representations $(\pi,V)$ such that $V$ is generated by $V^K$. 
Then, the functor $V\longmapsto V^K$ gives an equivalence between $\mathrm{Rep}(G)_K$ 
and the category of (unital) $\mathcal{H}(G,K)$-modules,
where $\mathcal{H}(G,K)$ denotes the Hecke algebra consisting of bi-$K$-invariant compactly supported functions
on $G$ (see \cite[Proposition VI.10.6 (ii)]{MR2567785}).
Further, there exists a finite subset $\mathcal{B}(G)_K$ of $\mathcal{B}(G)$
such that $\mathrm{Rep}(G)_K=\prod_{\mathfrak{s}\in\mathcal{B}(G)_K}\mathrm{Rep}(G)_{\mathfrak{s}}$ 
(see \cite[Proposition VI.10.6 (i)]{MR2567785}).

Let us begin a proof of Theorem \ref{thm:coh-fin-length}. We put $G=\GSp_{2n}(\Q_p)$, $J=\GU(n,D)$
and $V=H^i_c(M_\infty)$ for simplicity.
First we take a nice compact open subgroup $H$ of $J$ such that $\rho^{H}\neq 0$.
Then we have a decomposition $\mathrm{Rep}(J)=\mathrm{Rep}(J)_{H}\times \prod_{\mathfrak{s}\in\mathcal{B}(J)\setminus \mathcal{B}(J)_{H}}\mathrm{Rep}(G)_{\mathfrak{s}}$. Let $V=V_0\oplus V_1$ be the corresponding decomposition. Clearly the $J$-subrepresentation $V_0\subset V$
is also stable under the action of $G$.
As $\rho\in \mathrm{Rep}(J)_{H}$, we have
$\Ext^r_J(V,\rho)=\Ext^r_J(V_0,\rho)=\Ext^r_{\mathcal{H}(J,H)}(V_0^{H},\rho^{H})$.
By Corollary \ref{cor:finiteness+duality}, $V^H$ is a finitely generated $G$-module.
Therefore, the direct factor $V_0^{H}$ is also a finitely generated $G$-module.
Therefore we may take a nice compact open subgroup $K$ of $G$ such that $V_0^{H}\in \mathrm{Rep}(G)_K$.
We need the following lemma:

\begin{lem}\label{lem:Ext-Bernstein}
 Under the notation as above, $\Ext^r_J(V,\rho)^{\text{$\mathcal{D}_c$-$\mathrm{sm}$}}$ lies in the category $\mathrm{Rep}(G)_K$.
\end{lem}

\begin{prf}
 For simplicity, we put $R=\mathcal{H}(J,H)$. 
 The vector space $V_0^H$ has commuting actions of $G$ and $R$. We will call such structure a $(G,R)$-module.
 First we will prove that there exists a surjection $P\twoheadlongrightarrow V_0^H$ of $(G,R)$-modules
 from a $(G,R)$-module $P$ which lies in $\mathrm{Rep}(G)_K$ as a $G$-module and is projective as an $R$-module.
 Take an arbitrary surjection of $R$-modules $\varphi\colon Q\twoheadlongrightarrow V_0^H$
 from a free $R$-module $Q$.
 Put $Q'=\mathcal{H}(G)\otimes_\C Q$, where $\mathcal{H}(G)$ denotes the set of locally constant compactly supported functions on $G$.
 It has a structure of a $(G,R)$-module; the $G$-action is given by the left translation on $\mathcal{H}(G)$.
 Moreover, $Q'$ lies in $\mathrm{Rep}(G)$ as a $G$-module and is projective as an $R$-module.
 Define a map $Q'\longrightarrow V_0^H$ by $f\otimes x\longmapsto \int_G f(g)g\varphi(x)dg$.
 It is easy to see that this map is a surjection of $(G,R)$-modules.
 Let $P$ be the $\mathrm{Rep}(G)_K$-part of $Q'$. 
 Obviously $P$ is a direct factor of $Q'$ in the category of $(G,R)$-modules,
 and therefore it is projective as an $R$-module.
 The map $P\longrightarrow V_0^H$ induced from $Q'\longrightarrow V_0^H$ gives a desired surjection.
 
 As the kernel of this surjection is again a $(G,R)$-module belonging to $\mathrm{Rep}(G)_K$ as a $G$-module,
 we can repeat this construction and obtain a resolution $P_\bullet\longrightarrow V_0^H$.
 Let $K'$ be an arbitrary compact open subgroup of $G$.
 Since the functor $W\longmapsto e_{K'}W$ from the category of $\mathcal{D}_c$-modules to the category of
 $\C$-vector spaces is exact, we have
 \begin{align*}
  e_{K'}\Ext^r_J(V,\rho)&=e_{K'}\Ext^r_R(V_0^{H},\rho^{H})=e_{K'}H^r(\Hom_R(P_\bullet,\rho^{H}))\\
  &=H^r(e_{K'}\Hom_R(P_{\bullet},\rho^{H}))=H^r(\Hom_R(P_{\bullet},\rho^{H})^{K'}).
 \end{align*}
 By taking the inductive limit, we have 
 $\Ext^r_J(V,\rho)^{\text{$\mathcal{D}_c$-$\mathrm{sm}$}}=H^r(\Hom_R(P_{\bullet},\rho^{H})^{\mathrm{sm}})$. 
 Therefore it suffices to show that $\Hom_R(P_j,\rho^{H})^{\mathrm{sm}}$ lies in $\mathrm{Rep}(G)_K$ for each $j$.
 Since $\rho$ is an admissible representation of $J$, as $G$-representations we have
 \[
  \Hom_R(P_j,\rho^{H})^{\mathrm{sm}}\subset \Hom_{\C}(P_j,\rho^{H})^{\mathrm{sm}}\cong (P_j^\vee)^{\dim_\C\rho^{H}},
 \]
 where $(-)^\vee$ denotes the contragredient. By \cite[Proposition VI.10.6 (iii)]{MR2567785}, $P_j^\vee$ lies in $\mathrm{Rep}(G)_K$,
 hence $\Hom_R(P_j,\rho^{H})^{\mathrm{sm}}$ also lies in $\mathrm{Rep}(G)_K$. This concludes the proof.
\end{prf}

Now we can give a proof of Theorem \ref{thm:coh-fin-length}.
By a similar method as in the proof of \cite[Lemma 3.1]{LJLC}, we can show that 
\[
 (\Ext^r_J(V,\rho)^{\text{$\mathcal{D}_c$-$\mathrm{sm}$}})^{K'}=e_{K'}\Ext^r_J(V,\rho)\cong\Ext^r_J(V^{K'},\rho)
\]
is finite-dimensional for every compact open subgroup $K'$ of $G$.
This implies that $\Ext^r_J(V,\rho)^{\text{$\mathcal{D}_c$-$\mathrm{sm}$}}$ is admissible.
On the other hand, by Lemma \ref{lem:Ext-Bernstein}, we know that $\Ext^r_J(V,\rho)^{\text{$\mathcal{D}_c$-$\mathrm{sm}$}}$ is generated by
$\Ext^r_J(V^K,\rho)$. 
In particular $\Ext^r_J(V,\rho)^{\text{$\mathcal{D}_c$-$\mathrm{sm}$}}$ is finitely generated.
Therefore, \cite[Th\'eor\`eme VI.6.3]{MR2567785} tells us that
the $G$-representation $\Ext^r_J(V,\rho)^{\text{$\mathcal{D}_c$-$\mathrm{sm}$}}$ has finite length.
This completes the proof.

\begin{rem}
 The same argument works well if one switches the roles of $\GSp_{2n}$ and $\GU(n,D)$.
 Therefore, for an irreducible smooth representation $\pi$ of $\GSp_{2n}(\Q_p)$,
 the $\GU(n,D)$-representation $\Ext^r_{\GSp_{2n}(\Q_p)}(H^i_c(M_\infty),\pi)^{\text{$\mathcal{D}_c$-sm}}$
 has finite length, where in this case $\mathcal{D}_c$ denotes the convolution algebra of distributions
 on $\GU(n,D)$ with compact support.
\end{rem}

\noindent{\bfseries Acknowledgment}\quad
This work was supported by JSPS KAKENHI Grant Numbers 24740019, 15H03605.

\def\cftil#1{\ifmmode\setbox7\hbox{$\accent"5E#1$}\else
  \setbox7\hbox{\accent"5E#1}\penalty 10000\relax\fi\raise 1\ht7
  \hbox{\lower1.15ex\hbox to 1\wd7{\hss\accent"7E\hss}}\penalty 10000
  \hskip-1\wd7\penalty 10000\box7}
  \def\cftil#1{\ifmmode\setbox7\hbox{$\accent"5E#1$}\else
  \setbox7\hbox{\accent"5E#1}\penalty 10000\relax\fi\raise 1\ht7
  \hbox{\lower1.15ex\hbox to 1\wd7{\hss\accent"7E\hss}}\penalty 10000
  \hskip-1\wd7\penalty 10000\box7}
  \def\cftil#1{\ifmmode\setbox7\hbox{$\accent"5E#1$}\else
  \setbox7\hbox{\accent"5E#1}\penalty 10000\relax\fi\raise 1\ht7
  \hbox{\lower1.15ex\hbox to 1\wd7{\hss\accent"7E\hss}}\penalty 10000
  \hskip-1\wd7\penalty 10000\box7}
  \def\cftil#1{\ifmmode\setbox7\hbox{$\accent"5E#1$}\else
  \setbox7\hbox{\accent"5E#1}\penalty 10000\relax\fi\raise 1\ht7
  \hbox{\lower1.15ex\hbox to 1\wd7{\hss\accent"7E\hss}}\penalty 10000
  \hskip-1\wd7\penalty 10000\box7} \def\cprime{$'$} \def\cprime{$'$}
  \newcommand{\dummy}[1]{}
\providecommand{\bysame}{\leavevmode\hbox to3em{\hrulefill}\thinspace}
\providecommand{\MR}{\relax\ifhmode\unskip\space\fi MR }
\providecommand{\MRhref}[2]{%
  \href{http://www.ams.org/mathscinet-getitem?mr=#1}{#2}
}

\providecommand{\href}[2]{#2}

\end{document}